%Paper: funct-an/9411004
%From: piccoli@tsmi19.sissa.it
%Date: Thu, 03 Nov 1994 13:40:28 +0200

\magnification=\magstep1
\null
\hfuzz=20pt

\font\medbf=cmbx10 scaled \magstep2
\def \n{\noindent}
\def \v{\vskip 1em}
\def \vs{\vskip 2em}
\def\sgn{\hbox{sgn}}
\def \vsk{\vskip 4em}
\def \M {{\cal M}}

\def\DB{\Delta_B}
\def\DA{\Delta_A}
\def \R {I\!\!R}
\def\ve{\varepsilon}
\def \C {{\cal C}}

\def\H {{\cal H}}

%%%%amssym.def
%% @texfile{
%%     filename="amssym.def",
%%     version="2.1",
%%     date="5-APR-1991",
%%     filetype="TeX: option",
%%     copyright="Copyright (C) American Mathematical Society,
%%            all rights reserved.  Copying of this file is
%%            authorized only if either:
%%            (1) you make absolutely no changes to your copy
%%                including name; OR
%%            (2) if you do make changes, you first rename it to some
%%                other name.",
%%     author="American Mathematical Society",
%%     address="American Mathematical Society,
%%            Technical Support Group,
%%            P. O. Box 6248,
%%            Providence, RI 02940,
%%            USA",
%%     telephone="401-455-4080 or (in the USA) 800-321-4AMS",
%%     email="Internet: Tech-Support@Math.AMS.com",
%%     codetable="ISO/ASCII",
%%     checksumtype="line count",
%%     checksum="108",
%%     keywords="amsfonts, tex",
%%     abstract="This file contains definitions that perform the same
%%            functions as similar ones in AMS-TeX, so that the file
%%            AMSSYM.TEX can be used outside of AMS-TeX. Instructions
%%            for using this file and the AMS symbol fonts are
%%            included in the AMSFonts 2.0 User's Guide."
%%     }
%
%%%%%%%%%%%%%%%%%%%%%%%%%%%%%%%%%%%%%%%%%%%%%%%%%%%%%%%%%%%%%%%%%%%%%%%%
\expandafter\ifx\csname amssym.def\endcsname\relax \else\endinput\fi
%
%  Store the catcode of the @ in the csname so that it can be restored later.
\expandafter\edef\csname amssym.def\endcsname{%
       \catcode`\noexpand\@=\the\catcode`\@\space}
%  Set the catcode to 11 for use in private control sequence names.
\catcode`\@=11
%
%  Include all definitions related to the fonts msam, msbm and eufm, so that
%  when this file is used by itself, the results with respect to those fonts
%  are equivalent to what they would have been using AMS-TeX.
%  Most symbols in fonts msam and msbm are defined using \newsymbol;
%  however, a few symbols that replace composites defined in plain must be
%  defined with \mathchardef.

\def\undefine#1{\let#1\undefined}
\def\newsymbol#1#2#3#4#5{\let\next@\relax
 \ifnum#2=\@ne\let\next@\msafam@\else
 \ifnum#2=\tw@\let\next@\msbfam@\fi\fi
 \mathchardef#1="#3\next@#4#5}
\def\mathhexbox@#1#2#3{\relax
 \ifmmode\mathpalette{}{\m@th\mathchar"#1#2#3}%
 \else\leavevmode\hbox{$\m@th\mathchar"#1#2#3$}\fi}
\def\hexnumber@#1{\ifcase#1 0\or 1\or 2\or 3\or 4\or 5\or 6\or 7\or 8\or
 9\or A\or B\or C\or D\or E\or F\fi}

\font\tenmsa=msam10
\font\sevenmsa=msam7
\font\fivemsa=msam5
\newfam\msafam
\textfont\msafam=\tenmsa
\scriptfont\msafam=\sevenmsa
\scriptscriptfont\msafam=\fivemsa
\edef\msafam@{\hexnumber@\msafam}
\mathchardef\dabar@"0\msafam@39
\def\dashrightarrow{\mathrel{\dabar@\dabar@\mathchar"0\msafam@4B}}
\def\dashleftarrow{\mathrel{\mathchar"0\msafam@4C\dabar@\dabar@}}

\def\ulcorner{\delimiter"4\msafam@70\msafam@70 }
\def\urcorner{\delimiter"5\msafam@71\msafam@71 }
\def\llcorner{\delimiter"4\msafam@78\msafam@78 }
\def\lrcorner{\delimiter"5\msafam@79\msafam@79 }
\def\yen{{\mathhexbox@\msafam@55 }}
\def\checkmark{{\mathhexbox@\msafam@58 }}
\def\circledR{{\mathhexbox@\msafam@72 }}
\def\maltese{{\mathhexbox@\msafam@7A }}

\font\tenmsb=msbm10
\font\sevenmsb=msbm7
\font\fivemsb=msbm5
\newfam\msbfam
\textfont\msbfam=\tenmsb
\scriptfont\msbfam=\sevenmsb
\scriptscriptfont\msbfam=\fivemsb
\edef\msbfam@{\hexnumber@\msbfam}

\def\widehat#1{\setbox\z@\hbox{$\m@th#1$}%
 \ifdim\wd\z@>\tw@ em\mathaccent"0\msbfam@5B{#1}%
 \else\mathaccent"0362{#1}\fi}
\def\widetilde#1{\setbox\z@\hbox{$\m@th#1$}%
 \ifdim\wd\z@>\tw@ em\mathaccent"0\msbfam@5D{#1}%
 \else\mathaccent"0365{#1}\fi}
\font\teneufm=eufm10
\font\seveneufm=eufm7
\font\fiveeufm=eufm5
\newfam\eufmfam
\textfont\eufmfam=\teneufm
\scriptfont\eufmfam=\seveneufm
\scriptscriptfont\eufmfam=\fiveeufm

%  Restore the catcode value for @ that was previously saved.
\csname amssym.def\endcsname

\expandafter\ifx\csname pre amssym.tex at\endcsname\relax \else \endinput\fi
%  Otherwise we store the catcode of the @ in the csname.
\expandafter\chardef\csname pre amssym.tex at\endcsname=\the\catcode`\@
%  Set the catcode to 11 for use in private control sequence names.
\catcode`\@=11
%  Most symbols in fonts msam and msbm are defined using \newsymbol.  A few
%  that are delimiters or otherwise require special treatment have already
%  been defined as soon as the fonts were loaded.  Finally, a few symbols
%  that replace composites defined in plain must be undefined first.
\newsymbol\boxdot 1200
\newsymbol\boxplus 1201
\newsymbol\boxtimes 1202
\newsymbol\square 1003
\newsymbol\blacksquare 1004
\newsymbol\centerdot 1205
\newsymbol\lozenge 1006
\newsymbol\blacklozenge 1007
\newsymbol\circlearrowright 1308
\newsymbol\circlearrowleft 1309
\undefine\rightleftharpoons
\newsymbol\rightleftharpoons 130A
\newsymbol\leftrightharpoons 130B
\newsymbol\boxminus 120C
\newsymbol\Vdash 130D
\newsymbol\Vvdash 130E
\newsymbol\vDash 130F
\newsymbol\twoheadrightarrow 1310
\newsymbol\twoheadleftarrow 1311
\newsymbol\leftleftarrows 1312
\newsymbol\rightrightarrows 1313
\newsymbol\upuparrows 1314
\newsymbol\downdownarrows 1315
\newsymbol\upharpoonright 1316
 \let\restriction\upharpoonright
\newsymbol\downharpoonright 1317
\newsymbol\upharpoonleft 1318
\newsymbol\downharpoonleft 1319
\newsymbol\rightarrowtail 131A
\newsymbol\leftarrowtail 131B
\newsymbol\leftrightarrows 131C
\newsymbol\rightleftarrows 131D
\newsymbol\Lsh 131E
\newsymbol\Rsh 131F
\newsymbol\rightsquigarrow 1320
\newsymbol\leftrightsquigarrow 1321
\newsymbol\looparrowleft 1322
\newsymbol\looparrowright 1323
\newsymbol\circeq 1324
\newsymbol\succsim 1325
\newsymbol\gtrsim 1326
\newsymbol\gtrapprox 1327
\newsymbol\multimap 1328
\newsymbol\therefore 1329
\newsymbol\because 132A
\newsymbol\doteqdot 132B
 
\newsymbol\triangleq 132C
\newsymbol\precsim 132D
\newsymbol\lesssim 132E
\newsymbol\lessapprox 132F
\newsymbol\eqslantless 1330
\newsymbol\eqslantgtr 1331
\newsymbol\curlyeqprec 1332
\newsymbol\curlyeqsucc 1333
\newsymbol\preccurlyeq 1334
\newsymbol\leqq 1335
\newsymbol\leqslant 1336
\newsymbol\lessgtr 1337
\newsymbol\backprime 1038
\newsymbol\risingdotseq 133A
\newsymbol\fallingdotseq 133B
\newsymbol\succcurlyeq 133C
\newsymbol\geqq 133D
\newsymbol\geqslant 133E
\newsymbol\gtrless 133F
\newsymbol\sqsubset 1340
\newsymbol\sqsupset 1341
\newsymbol\vartriangleright 1342
\newsymbol\vartriangleleft 1343
\newsymbol\trianglerighteq 1344
\newsymbol\trianglelefteq 1345
\newsymbol\bigstar 1046
\newsymbol\between 1347
\newsymbol\blacktriangledown 1048
\newsymbol\blacktriangleright 1349
\newsymbol\blacktriangleleft 134A
\newsymbol\vartriangle 134D
\newsymbol\blacktriangle 104E
\newsymbol\triangledown 104F
\newsymbol\eqcirc 1350
\newsymbol\lesseqgtr 1351
\newsymbol\gtreqless 1352
\newsymbol\lesseqqgtr 1353
\newsymbol\gtreqqless 1354
\newsymbol\Rrightarrow 1356
\newsymbol\Lleftarrow 1357
\newsymbol\veebar 1259
\newsymbol\barwedge 125A
\newsymbol\doublebarwedge 125B
\undefine\angle
\newsymbol\angle 105C
\newsymbol\measuredangle 105D
\newsymbol\sphericalangle 105E
\newsymbol\varpropto 135F
\newsymbol\smallsmile 1360
\newsymbol\smallfrown 1361
\newsymbol\Subset 1362
\newsymbol\Supset 1363
\newsymbol\Cup 1264
 
\newsymbol\Cap 1265
 
\newsymbol\curlywedge 1266
\newsymbol\curlyvee 1267
\newsymbol\leftthreetimes 1268
\newsymbol\rightthreetimes 1269
\newsymbol\subseteqq 136A
\newsymbol\supseteqq 136B
\newsymbol\bumpeq 136C
\newsymbol\Bumpeq 136D
\newsymbol\lll 136E
 
\newsymbol\ggg 136F
 
\newsymbol\circledS 1073
\newsymbol\pitchfork 1374
\newsymbol\dotplus 1275
\newsymbol\backsim 1376
\newsymbol\backsimeq 1377
\newsymbol\complement 107B
\newsymbol\intercal 127C
\newsymbol\circledcirc 127D
\newsymbol\circledast 127E
\newsymbol\circleddash 127F
\newsymbol\lvertneqq 2300
\newsymbol\gvertneqq 2301
\newsymbol\nleq 2302
\newsymbol\ngeq 2303
\newsymbol\nless 2304
\newsymbol\ngtr 2305
\newsymbol\nprec 2306
\newsymbol\nsucc 2307
\newsymbol\lneqq 2308
\newsymbol\gneqq 2309
\newsymbol\nleqslant 230A
\newsymbol\ngeqslant 230B
\newsymbol\lneq 230C
\newsymbol\gneq 230D
\newsymbol\npreceq 230E
\newsymbol\nsucceq 230F
\newsymbol\precnsim 2310
\newsymbol\succnsim 2311
\newsymbol\lnsim 2312
\newsymbol\gnsim 2313
\newsymbol\nleqq 2314
\newsymbol\ngeqq 2315
\newsymbol\precneqq 2316
\newsymbol\succneqq 2317
\newsymbol\precnapprox 2318
\newsymbol\succnapprox 2319
\newsymbol\lnapprox 231A
\newsymbol\gnapprox 231B
\newsymbol\nsim 231C
\newsymbol\ncong 231D
\newsymbol\diagup 231E
\newsymbol\diagdown 231F
\newsymbol\varsubsetneq 2320
\newsymbol\varsupsetneq 2321
\newsymbol\nsubseteqq 2322
\newsymbol\nsupseteqq 2323
\newsymbol\subsetneqq 2324
\newsymbol\supsetneqq 2325
\newsymbol\varsubsetneqq 2326
\newsymbol\varsupsetneqq 2327
\newsymbol\subsetneq 2328
\newsymbol\supsetneq 2329
\newsymbol\nsubseteq 232A
\newsymbol\nsupseteq 232B
\newsymbol\nparallel 232C
\newsymbol\nmid 232D
\newsymbol\nshortmid 232E
\newsymbol\nshortparallel 232F
\newsymbol\nvdash 2330
\newsymbol\nVdash 2331
\newsymbol\nvDash 2332
\newsymbol\nVDash 2333
\newsymbol\ntrianglerighteq 2334
\newsymbol\ntrianglelefteq 2335
\newsymbol\ntriangleleft 2336
\newsymbol\ntriangleright 2337
\newsymbol\nleftarrow 2338
\newsymbol\nrightarrow 2339
\newsymbol\nLeftarrow 233A
\newsymbol\nRightarrow 233B
\newsymbol\nLeftrightarrow 233C
\newsymbol\nleftrightarrow 233D
\newsymbol\divideontimes 223E
\newsymbol\varnothing 203F
\newsymbol\nexists 2040
\newsymbol\Finv 2060
\newsymbol\Game 2061
\newsymbol\mho 2066
\newsymbol\eth 2067
\newsymbol\eqsim 2368
\newsymbol\beth 2069
\newsymbol\gimel 206A
\newsymbol\daleth 206B
\newsymbol\lessdot 236C
\newsymbol\gtrdot 236D
\newsymbol\ltimes 226E
\newsymbol\rtimes 226F
\newsymbol\shortmid 2370
\newsymbol\shortparallel 2371
\newsymbol\smallsetminus 2272
\newsymbol\thicksim 2373
\newsymbol\thickapprox 2374
\newsymbol\approxeq 2375
\newsymbol\succapprox 2376
\newsymbol\precapprox 2377
\newsymbol\curvearrowleft 2378
\newsymbol\curvearrowright 2379
\newsymbol\digamma 207A
\newsymbol\varkappa 207B
\newsymbol\Bbbk 207C
\newsymbol\hslash 207D
\undefine\hbar
\newsymbol\hbar 207E
\newsymbol\backepsilon 237F
%  Restore the catcode value for @ that was previously saved.
\catcode`\@=\csname pre amssym.tex at\endcsname

%\endinput
%%%%%%%%

%%%%fonts
%FONTS
\font\twelverm=cmr12
\font\twelvei=cmmi12
\font\twelvesy=cmsy10
\font\twelvebf=cmbx12
\font\twelvett=cmtt12
\font\twelveit=cmti12
\font\twelvesl=cmsl12

\font\ninerm=cmr9
\font\ninei=cmmi9
\font\ninesy=cmsy9
\font\ninebf=cmbx9
\font\ninett=cmtt9
\font\nineit=cmti9
\font\ninesl=cmsl9

\font\ninesc = cmcsc10 at 9pt

\font\eightrm=cmr8
\font\eighti=cmmi8
\font\eightsy=cmsy8
\font\eightbf=cmbx8
\font\eighttt=cmtt8
\font\eightit=cmti8
\font\eightsl=cmsl8

\font\sixrm=cmr6
\font\sixi=cmmi6
\font\sixsy=cmsy6
\font\sixbf=cmbx6

\catcode`@=11 % we will access private macros of plain TeX (carefully)
\newskip\ttglue

%MACRO TWELVEPOINT
\def\twelvepoint{\def\rm{\fam0\twelverm}% switch to 12-point type
\textfont0=\twelverm  \scriptfont0=\ninerm
\scriptscriptfont0=\sevenrm
\textfont1=\twelvei  \scriptfont1=\ninei  \scriptscriptfont1=\seveni
\textfont2=\twelvesy  \scriptfont2=\ninesy
\scriptscriptfont2=\sevensy
\textfont3=\tenex  \scriptfont3=\tenex  \scriptscriptfont3=\tenex
\textfont\itfam=\twelveit  \def\it{\fam\itfam\twelveit}%
\textfont\slfam=\twelvesl  \def\sl{\fam\slfam\twelvesl}%
\textfont\ttfam=\twelvett  \def\tt{\fam\ttfam\twelvett}%
\textfont\bffam=\twelvebf  \scriptfont\bffam=\ninebf
\scriptscriptfont\bffam=\sevenbf  \def\bf{\fam\bffam\twelvebf}%
\tt  \ttglue=.5em plus.25em minus.15em
\normalbaselineskip=15pt
\setbox\strutbox=\hbox{\vrule height10pt depth5pt width0pt}%
\let\sc=\tenrm  \let\big=\twelvebig  \normalbaselines\rm}

%MACRO TENPOINT
\def\tenpoint{\def\rm{\fam0\tenrm}% switch to 10-point type
\textfont0=\tenrm  \scriptfont0=\sevenrm  \scriptscriptfont0=\fiverm
\textfont1=\teni  \scriptfont1=\seveni  \scriptscriptfont1=\fivei
\textfont2=\tensy  \scriptfont2=\sevensy  \scriptscriptfont2=\fivesy
\textfont3=\tenex  \scriptfont3=\tenex  \scriptscriptfont3=\tenex
\textfont\itfam=\tenit  \def\it{\fam\itfam\tenit}%
\textfont\slfam=\tensl  \def\sl{\fam\slfam\tensl}%
\textfont\ttfam=\tentt  \def\tt{\fam\ttfam\tentt}%
\textfont\bffam=\tenbf  \scriptfont\bffam=\sevenbf
\scriptscriptfont\bffam=\fivebf  \def\bf{\fam\bffam\tenbf}%
\tt  \ttglue=.5em plus.25em minus.15em
\normalbaselineskip=12pt
\setbox\strutbox=\hbox{\vrule height8.5pt depth3.5pt width0pt}%
\let\sc=\eightrm  \let\big=\tenbig  \normalbaselines\rm}

%MACRO NINEPOINT
\def\ninepoint{\def\rm{\fam0\ninerm}% switch to 9-point type
\textfont0=\ninerm  \scriptfont0=\sixrm  \scriptscriptfont0=\fiverm
\textfont1=\ninei  \scriptfont1=\sixi  \scriptscriptfont1=\fivei
\textfont2=\ninesy  \scriptfont2=\sixsy  \scriptscriptfont2=\fivesy
\textfont3=\tenex  \scriptfont3=\tenex  \scriptscriptfont3=\tenex
\textfont\itfam=\nineit  \def\it{\fam\itfam\nineit}%
\textfont\slfam=\ninesl  \def\sl{\fam\slfam\ninesl}%
\textfont\ttfam=\ninett  \def\tt{\fam\ttfam\ninett}%
\textfont\bffam=\ninebf  \scriptfont\bffam=\sixbf
\scriptscriptfont\bffam=\fivebf  \def\bf{\fam\bffam\ninebf}%
\tt  \ttglue=.5em plus.25em minus.15em
\normalbaselineskip=11pt
\setbox\strutbox=\hbox{\vrule height8pt depth3pt width0pt}%
\let\sc=\sevenrm  \let\big=\ninebig  \normalbaselines\rm}

%MACRO EIGHTPOINT
\def\eightpoint{\def\rm{\fam0\eightrm}% switch to 8-point type
\textfont0=\eightrm  \scriptfont0=\sixrm  \scriptscriptfont0=\fiverm
\textfont1=\eighti  \scriptfont1=\sixi  \scriptscriptfont1=\fivei
\textfont2=\eightsy  \scriptfont2=\sixsy  \scriptscriptfont2=\fivesy
\textfont3=\tenex  \scriptfont3=\tenex  \scriptscriptfont3=\tenex
\textfont\itfam=\eightit  \def\it{\fam\itfam\eightit}%
\textfont\slfam=\eightsl  \def\sl{\fam\slfam\eightsl}%
\textfont\ttfam=\eighttt  \def\tt{\fam\ttfam\eighttt}%
\textfont\bffam=\eightbf  \scriptfont\bffam=\sixbf
\scriptscriptfont\bffam=\fivebf  \def\bf{\fam\bffam\eightbf}%
\tt  \ttglue=.5em plus.25em minus.15em
\normalbaselineskip=9pt
\setbox\strutbox=\hbox{\vrule height7pt depth2pt width0pt}%
\let\sc=\sixrm  \let\big=\eightbig  \normalbaselines\rm}

%MACRO BIG
\def\twelvebig#1{{\hbox{$\textfont0=\twelverm\textfont2=\twelvesy
	\left#1\vbox to10pt{}\right.\n@space$}}}
\def\tenbig#1{{\hbox{$\left#1\vbox to8.5pt{}\right.\n@space$}}}
\def\ninebig#1{{\hbox{$\textfont0=\tenrm\textfont2=\tensy
	\left#1\vbox to7.25pt{}\right.\n@space$}}}
\def\eightbig#1{{\hbox{$\textfont0=\ninerm\textfont2=\ninesy
	\left#1\vbox to6.5pt{}\right.\n@space$}}}

%for 10-point math in 9-point territory
%%%%macro

%
\nopagenumbers
\headline={\ifodd\pageno\rightheadline \else\leftheadline\fi}
\def\rightheadline{\hfil{\eightpoint\titolo}
\hfil\tenrm\folio}
\def\leftheadline{\tenrm\folio\hfil{\eightpoint\autori}
\hfil} \topskip=25pt

\newcount\parn \newcount\forn \newcount\teon
%\newcount\lemn \newcount\pron
%\newcount\esen \newcount\ossn \newcount\defn
\newcount\vol \newcount\pag
%
%\font\mate=msym10
%
\parn=0
\def\newpar#1{{\global\forn=0 \global\teon=0
\global\advance\parn by 1} \bigskip {\noindent\bf \the\parn
. #1} \bigskip}  %
\def\phor{{\global\advance\forn by 1}
\leqno{(\the\parn.\the\forn)}}
\long\def\theo#1{{\global\advance\teon by 1}\medskip {\noindent\bf Theorem
\the\parn.\the\teon}\quad{\sl #1}\medskip}
\long\def\prop#1{{\global\advance\teon by 1}\medskip {\noindent\bf
Proposition \the\parn.\the\teon}\quad{\sl #1}\medskip}
\long\def\lemma#1{{\global\advance\teon by 1}\medskip {\noindent\bf
Lemma \the\parn.\the\teon}\quad{\sl #1}\medskip}
\def\ehse{{\global\advance\teon by 1}\medskip {\noindent\bf Example
\the\parn.\the\teon}\quad}
\def\ohss{{\global\advance\teon by1}\medskip {\noindent\bf Remark
\the\parn.\the\teon}\quad}
\def\dhef{{\global\advance\teon by 1}\medskip {\noindent\bf
Definition \the\parn.\the\teon}\quad}
\def\sqr#1#2{{\vcenter{\vbox{\hrule height.#2pt
													\hbox{\vrule width.#2pt height#1pt \kern#1pt
													\vrule width.#2pt}
             \hrule height.#2pt}}}}

\def\norma#1#2{\|#1\|_{\lower 4pt \hbox{$\scriptstyle #2$}}}

\def\bibart#1#2#3#4#5#6#7#8{\global\vol=#5 \global\pag=#7 \par
\hangindent=20pt \hangafter=-10 {\item{[{\bf #1}]}}{\ninesc #2}: {#3},
{\sl #4},%
{\ifnum\vol=0\else{{\bf\ #5},}\fi}
{#6}%
{\ifnum\pag=0\else{, p.~{#7}--{#8}}\fi}%
. \medskip}
\mathchardef\emptyset="001F

\nopagenumbers
\headline={\ifodd\pageno\rightheadline \else\leftheadline\fi}
\def\rightheadline{\hfil{\eightpoint\titolo}
\hfil\tenrm\folio}
\def\leftheadline{\tenrm\folio\hfil{\eightpoint\autori}
\hfil} \topskip=25pt

\newcount\parn \newcount\forn \newcount\teon
%\newcount\lemn \newcount\pron
%\newcount\esen \newcount\ossn \newcount\defn
\newcount\vol \newcount\pag
%
%\font\mate=msym10
%
\parn=0
\def\newpar#1{{\global\forn=0 \global\teon=0
\global\advance\parn by 1} \bigskip {\noindent\bf \the\parn
. #1} \bigskip}  %
\def\phor{{\global\advance\forn by 1}
\leqno{(\the\parn.\the\forn)}}
\long\def\theo#1{{\global\advance\teon by 1}\medskip {\noindent\bf Theorem
\the\parn.\the\teon}\quad{\sl #1}\medskip}
\long\def\prop#1{{\global\advance\teon by 1}\medskip {\noindent\bf
Proposition \the\parn.\the\teon}\quad{\sl #1}\medskip}
\long\def\lemma#1{{\global\advance\teon by 1}\medskip {\noindent\bf
Lemma \the\parn.\the\teon}\quad{\sl #1}\medskip}
\def\ehse{{\global\advance\teon by 1}\medskip {\noindent\bf Example
\the\parn.\the\teon}\quad}
\def\ohss{{\global\advance\teon by1}\medskip {\noindent\bf Remark
\the\parn.\the\teon}\quad}
\def\dhef{{\global\advance\teon by 1}\medskip {\noindent\bf
Definition \the\parn.\the\teon}\quad}
\def\sqr#1#2{{\vcenter{\vbox{\hrule height.#2pt
													\hbox{\vrule width.#2pt height#1pt \kern#1pt
													\vrule width.#2pt}
             \hrule height.#2pt}}}}

\def\norma#1#2{\|#1\|_{\lower 4pt \hbox{$\scriptstyle #2$}}}

\def\bibart#1#2#3#4#5#6#7#8{\global\vol=#5 \global\pag=#7 \par
\hangindent=20pt \hangafter=-10 {\item{[{\bf #1}]}}{\ninesc #2}: {#3},
{\sl #4},%
{\ifnum\vol=0\else{{\bf\ #5},}\fi}
{#6}%
{\ifnum\pag=0\else{, p.~{#7}--{#8}}\fi}%
. \medskip}
\mathchardef\emptyset="001F
%%%%

\tenpoint
\def\autori{B. PICCOLI}
\def\titolo{Time-Optimal Control Problems for the Swing and the Ski}
\baselineskip = 18 truept

\n {\bf 1.Introduction.}
\vs
This paper is concerned with a class of time-optimal control problems
for the swing and the ski.

\n We first consider the motion of a man standing on a swing.
For simplicity, we neglect friction and air resistance and
assume that the mass of the swinger is concentrated in his baricenter $B$.
Let $\theta$ be the angle formed by the swing and the downward
vertical direction and $r$ the radius of oscillation, i.e. the
distance between the baricenter $B$ and the center of rotation $O$.
Assuming that the mass is normalized to a unit, the corresponding
Lagrangean function takes the form:
$$L\;=\;{1\over 2}\;(\dot r^2+r^2\;\dot\theta^2)+g\;r\;cos(\theta)
\eqno(1.1)$$
where $g$ is the gravity acceleration.

\n We assume that, by bending his knees, the swinger can vary his
radius of oscillation.
This amounts to the addition of a constraint $r=u(t)$, implemented
by forces acting on $B$, parallel to the vector $OB$. The function
$u(t)$ can be regarded here as a control, whose values are chosen
by the man riding on the swing, within certain physical bounds,
say $u(t)\in[r_-,r_+]$ with $0<r_-<r_+$.

\n Writing (1.1) as a first order
system for $(\theta,\dot\theta)$ we obtain an impulsive control system.
On the other hand, using the coordinates $x_1=\theta,\;
x_2=\dot\theta r^2$ (the angular momentum), following Alberto
Bressan (1993), we obtain the
nonimpulsive control system:
$$\cases{\dot x_1\;=\;x_2/v^2 &\cr
                               & $v\in[r_-,r_+]$\cr
          \dot x_2\;=\;-g\;v\;sin(x_1)              &\cr}.\eqno(1.2)$$

\n For the ski, we use the model considered in Aldo Bressan (1991), Bressan
and Motta (1994), with a special approximation of the skier.

Our main concern is the existence and the structure of time-optimal
controls. As a first step, we introduce an auxiliary control system $\Sigma$
with control entering linearly and establish the relationships between
this and the original system.
We then perform a detailed study of the auxiliary system, using
the geometric techniques developed in Sussmann (1987 a,b,c)
and in Piccoli (1993 a,b).
These geometric techniques permit us to know the structure of the time
optimal trajectories in generic cases and possibly to solve
explicitely some optimization problem. See section 4 for examples.
In turn this provides accurate information on the time-optimal controls
for the original system (1.2).

Since the set of admissible velocities for $\Sigma$
is always convex, a time optimal
control exists for general boundary conditions. However, we prove
that a time optimal control for (1.2) does exist only when
the corresponding optimal control
for $\Sigma$ is bang-bang. In this case the optimal trajectory is the same
for the two control systems. We show that, for every
control constraint $v\in [r_-,r_+]$ with $0<r_-<r_+$, there
exists some pair of points $x,\tilde x\in\R^2$ such that
no control for (1.2) steers $x$ to $\tilde x$ in minimum time.

\n Then we consider the problem of reaching with the swing
a given angle $\bar\theta$ in minimum time,
having assigned initial condition $\theta,\dot\theta\;r^2$.
This can be stated as a Mayer problem for (1.2).
We show how to solve this problem under the assumptions
$|\theta|,|\bar\theta|\leq \pi/2$.
To solve numerically this problem, it sufficies to find a suitable class
of solutions to (1.2) with constant control $v\in \{r_-,r_+\}$,
and then to solve the complementary linear system (2.6) for each trajectory.
The optimal trajectory is carachterized by the (final) transversality
condition.
We also show that, for a special
class of boundary conditions and data $r_{\pm}$, no optimal solution exists.

\n For the ski model we obtain similar results under the assumption
that the curvature of the ski trail is constant.

\n For an introduction to impulsive control systems we refer to
Alberto Bressan (1993).
\v
{\bf Acknoledgments.} We wish to thank Prof. Alberto Bressan and
Prof. Aldo Bressan for having suggested the problem and for many
useful advises.
\vsk
\n {\bf 2.Basic definitions.}
\vs
By a {\it curve} in $\R^n$ we mean a continuous map $\gamma:I\mapsto\R^n$,
where $I$ is some real interval. We use the symbol $Dom$ for the
domain so that if $\gamma:I\mapsto\R^n$ then $Dom(\gamma)=I$.
We use the symbol $\gamma\restriction J$, where
$J\subset Dom(\gamma)$ is an interval, to denote the restriction of $\gamma$
to $J$.

A ${\C}^1$ {\it vector field} on $\R^2$ is a continuosly differentiable
map $F:\R^2\rightarrow\R^2$. It can be written
in the form:
$$F\ =\ \alpha\ {\partial}_x\ +\ \beta\ {\partial}_y\eqno(2.1)$$
where ${\partial}_x,\ {\partial}_y$ are the constant vector fields
with components (1,0), (0,1) respectively.
If we have the representation (2.1) for a vector field $F$ then
we use the symbol $\nabla F$ to denote the $2\times 2$ Jacobian matrix:
$$\nabla F=\left(\matrix
{ {\partial\alpha\over\partial x} & {\partial\alpha\over\partial y}\cr
                               &  \cr
  {\partial\beta\over\partial x} & {\partial\beta\over\partial y}  \cr }
\right).$$
\n The {\it Lie-bracket} of two vector fields $F,\ G$ is the vector
field:
$$[F,G]\quad =\quad \nabla G\cdot F\quad-\quad\nabla F\cdot G.\eqno(2.2)$$

We consider a two dimensional autonomous control system:
$$\dot x\ =\ h(x,u)\qquad u\in U\eqno(2.3)$$
where $U\subset\R^m$ is compact and $h\in {\C}^1(\R^2\times\R^m,\R^2)$.
A {\it control} is a measurable function $u:[a,b]\mapsto U$
where $-\infty<a\leq b<+\infty$. As for the curves we use the symbol
$Dom$ for the domain.
A {\it trajectory } for a control $u$ is an absolutely
continuous curve $\gamma :Dom(u)\mapsto\R^2$ which satisfies the
equation:
$$\dot\gamma(t)\quad=\quad h(\gamma(t),u(t))$$
for almost every $t\in Dom(u)$.

If $\gamma:[a,b]\mapsto\R^2$ is a trajectory of (2.3) we use the symbol
$In(\gamma)$ to denote its initial point $\gamma(a)$
and $Term(\gamma)$ to denote its terminal point $\gamma(b)$. We define:
$$T(\gamma)=b-a.$$
i.e. $T(\gamma)$ is the time along $\gamma$.

A trajectory $\gamma$ is said to be {\it time optimal}
if for every trajectory $\gamma'$ with $In(\gamma')=
In(\gamma)$ and $Term(\gamma')=Term(\gamma)$, we have that
$T(\gamma')\geq T(\gamma)$.

If $u_1:[a,b]\mapsto U$ and $u_2:[b,c]\mapsto U$
are controls, we use $u_2\ast u_1$ to denote the control defined by:
$$(u_2*u_1)(t)=\cases{ u_1(t) &for $t\in Dom(u_1)$ \cr
                              &                    \cr
                       u_2(t) &for $t\in Dom(u_2)$ \cr}.$$
This control is called the {\it concatenation} of $u_1$ and $u_2$.

\n If $\gamma_1:[a,b]\mapsto\R^2,\ \gamma_2:[b,c]\mapsto\R^2$
are trajectories for $u_1$ and $u_2$ such that
$\gamma_1(b)=\gamma_2(b)$, then the {\it concatenation}
$\gamma_2*\gamma_1$ is the trajectory:
$$(\gamma_2*\gamma_1)(t)=\cases{\gamma_1(t) &for $t\in Dom(\gamma_1)$ \cr
                                            &                         \cr
                                \gamma_2(t) &for $t\in Dom(\gamma_2)$ \cr}.$$
Now we consider a control system with control entering linearly:
$$\dot x\ =\ F(x)+u\;G(x)\qquad |u|\leq 1\eqno(2.4)$$
where $F,G$ are two ${\C}^1$ vector field on $\R^2$ and $u$ is scalar.
We define:
$$X=F-G\qquad Y=F+G.$$
An $X-$trajectory is a trajectory
corresponding to a constant control $u$ whose value is equal
to $-1$. We define $Y-$trajectories in a similar way, using the control
$u=+1$ rather than $u=-1$.
An $X*Y-$trajectory
is concatenation of a $Y-$trajectory and an $X-$trajectory ( the
$Y-$trajectory comes first), and similarly is defined
a $Y*X-$trajectory.

\n A {\it bang-bang trajectory} is a trajectory that is a concatenation
of $X-$ and $Y-$trajectories.

\n A time $t\in Dom(\gamma)$ is called
a {\it switching time} for $\gamma$ if, for each $\varepsilon>0$,
$\gamma\restriction [t-\varepsilon,~t+\varepsilon]$
is neither an $X$-trajectory nor a $Y$-trajectory.
If $t$ is a switching time for $\gamma$ then
we say that $\gamma(t)$ is a {\it switching point for $\gamma$},
or that $\gamma$ has a switching at $\gamma(t)$.

%%%----------------PMP------------------

For the rest of this section we call $\Sigma$ the system (2.4).
An {\it admissible pair} for the system $\Sigma$ is a couple $(u,\gamma)$
such that $u$ is a control and $\gamma$ is a trajectory corresponding to $u$.
We use the symbol $Adm(\Sigma)$ to denote the set of admissible pairs
and we say that $(u,\gamma)\in Adm(\Sigma)$
is optimal if $\gamma$ is optimal.

A {\it variational vector field along $(u,\gamma)\in Adm(\Sigma)$}
is a vector-valued absolutely continuous function
$v:Dom(\gamma)\mapsto\R^2$ that satisfies the equation:
$$\dot v(t)=\Big((\nabla F)(\gamma(t))+u(t)(\nabla G)(\gamma (t))\Big)
\cdot v(t)\eqno(2.5)$$
for almost all $t\in Dom(\gamma)$.

A {\it variational covector field along $(u,\gamma)\in Adm(\Sigma)$}
is an absolutely continuous function
$\lambda:Dom(\gamma)\mapsto\R^2_*$ that satisfies the equation:
$$\dot \lambda(t)=-\lambda(t)\cdot\Big((\nabla F)(\gamma(t))+u(t)
(\nabla G)(\gamma (t))\Big)
\eqno(2.6)$$
for almost all $t\in Dom(\gamma)$. Here $\R^2_*$ denotes the space of
row vectors. In Sussmann (1987 a) it was proved:
\v
\n {\bf Lemma 2.1} {\it Let $(u,\gamma)\in Adm(\Sigma)$ and
let $\lambda:Dom(\gamma)\mapsto\R^2_*$ be absolutely continuous.
Then $\lambda$ is a variational covector field along $(u,\gamma)$
if and only if the function $t\mapsto\lambda(t)\cdot v(t)$
is constant for every variational vector field $v$ along $(u,\gamma)$.}
\v
The Hamiltonian $\H:\R^2_*\times\R^2\times\R\mapsto\R$ is defined as
$$\H(\lambda,x,u)=\lambda\cdot(F(x)+uG(x)).\eqno(2.7)$$
If $\lambda$ is a variational covector field along $(u,\gamma)\in Adm(\Sigma)$,
we say that $\lambda$ is {\it maximizing} if:
$$\H(\lambda(t),\gamma(t),u(t))=\max\ \{\H(\lambda(t),\gamma(t),w):
|w|\leq 1\}\eqno(2.8)$$
for almost all $t\in Dom(\gamma)$.

The {\it Pontryagin Maximum Principle} (PMP) states
that, if $(u,\gamma)\in Adm(\Sigma)$
is time optimal, then there exists:
\v
\item{(PMP1)} {\it A non trivial maximizing
variational covector field $\lambda$
along $(u,\gamma)$}
\v
\item{(PMP2)} {\it A constant $\lambda_0\leq 0$ such that :
$\H(\lambda(t),\gamma(t),u(t))+\lambda_0=0$}
for almost all $t\in Dom(\gamma)$.
\v
In this case $\lambda$ is called
an {\it adjoint covector field along $(u,\gamma)$} or
simply an {\it adjoint variable}, and we say that
$(\gamma,\lambda)$ satisfies the PMP.

If $\lambda$ is an adjoint covector field along $(u,\gamma)\in Adm(\Sigma)$,
the corresponding {\it switching function} is defined as:
$$\phi_{\lambda}(t)=\lambda(t)\cdot G(\gamma(t)).\eqno(2.9)$$
{}From the above definition it follows:
\v
\n {\bf Lemma 2.2}
{\it
Let $(u,\gamma)\in Adm(\Sigma)$ be optimal and let
$\lambda$ be an adjoint covector field along $(u,\gamma)$.  Then:

\item {a)} The switching function $\phi_{\lambda}$ is continuous.

\item {b)} If $\phi(t)>0$ for all $t$ in some interval $I$, then $u(t)\equiv 1$
for almost all $t\in I$ and $\gamma\restriction I$ is a $Y$-trajectory.

\item {c)} If $\phi(t)<0$ for all $t$ in some interval $I$, then $u(t)\equiv
-1$
for almost all $t\in I$ and $\gamma\restriction I$ is an $X$-trajectory.

}
\v
For each $x\in \R^2$, one can
form the $2\times 2$ matrices whose columns are the vectors $F,~G,$ or $[F,G]$.
As in Sussmann (1987 a), we shall use the following scalar
functions on $\R^2$:
$$\DA(x)\doteq det\big(F(x),G(x)\big)\eqno(2.10)$$
$$\DB(x)\doteq det\big(G(x),[F,G](x)\big)\eqno(2.11)$$
where $det$ stands for determinant.

\n Consider $(u,\gamma)\in Adm(\Sigma)$, $t_0\in Dom(\gamma)$ and $v_0\in\R^2$.
We write $v(v_0,t_0;t)$ to denote the value at time $t$ of the
variational vector field
along $(u,\gamma)$ satisfying (2.5) together with the boundary condition
$v(t_0)=v_0$. In Piccoli (1993 a) it was proved the following:
\v
\n{\bf Lemma 2.3.}  {\it Let $(u,\gamma)\in Adm(\Sigma)$, $t_0\in Dom(\gamma)$,
and $v_0\in\R^2$, $v_0\not=0$.
For every $t$ such that $G(\gamma(t))\not=0$, define the angle:
$$\alpha(t)=\arg\bigg( v_0,~~v\big( G(\gamma(t)),t;~t_0\big)
\bigg),\eqno(2.12)$$
Then, one has:}
$$\sgn\big(\dot \alpha(t)\big)=\sgn\Big(\DB\big(\gamma(t)\big)\Big).
\eqno(2.13)$$
\v

A point $x\in\R^2$ is called an {\it ordinary point} if
$$\DA(x)\cdot \DB(x)\not= 0.$$
On the set of ordinary points we define the scalar functions $f$, $g$
as the coefficients of the linear combination
$$[F,G](x)=f(x)F(x)+g(x)G(x).\eqno(2.14)$$
By direct calculations we have:
$$f=-{\DB\over\DA}.$$

A point $x$ at which $\DA(x)\DB(x)=0$ is called a {\it nonordinary point}.
A {\it nonordinary arc} is a $\C^1$ one-dimensional connected embedded
submanifold $S$ of $\R^2$, with the property that every $x\in S$
is a nonordinary point.
A nonordinary arc will be said {\it isolated}, and
will be called an INOA, if there exists a set
$\Omega$ satisfying the following conditions:
\v
\item {(C1)} $\Omega$ is an open connected subset of $\R^2$

\item {(C2)} $S$ is a relatively closed subset of $\Omega$

\item {(C3)} If $x\in \Omega-S$ then $x$ is an ordinary point

\item {(C4)} The set $\Omega-S$ has exactly two connected components.
\v
\n A {\it turnpike} is an isolated
nonordinary arc that satisfies the following
conditions:
\v
\item {(S1)} For each $x\in S$ the vectors $X(x)$ and $Y(x)$ are not
tangent to $S$ and point to opposite sides of $S$

\item {(S2)} For each $x\in S$ one has $\Delta_B(x)=0$ and
$\Delta_A(x)\not= 0$

\item {(S3)} Let $\Omega$ be an open set which satisfies (C1)-(C4) above.
If $\Omega_X$ and $\Omega_Y$ are the connected components
of $\Omega-S$ labelled in such a way that $X(x)$ points into $\Omega_X$
and $Y(x)$ points into $\Omega_Y$, then the function $f$ in (2.14) satisfies
$$f(x)>0\qquad \hbox{on}\quad \Omega_Y\qquad\qquad\qquad f(x)<0
\qquad \hbox{on}\quad \Omega_X.$$
\v
Next, consider a turnpike $S$ and a point $x_0\in S$. We wish to construct
a trajectory $\gamma$ of (2.4) such that
$\gamma(t_0)=x_0$ and $\gamma(t)\in S$ for each $t\in Dom(\gamma)\doteq
[t_0,t_1]$.
Clearly, one should have
$\DB(\gamma(t)) \equiv 0$ for all $t$.  Since
$\DB(\gamma(t_0))=0$, it suffices to verify that:
$${d\over dt}\DB\big(\gamma(t)\big)=
(\nabla\DB\cdot\dot\gamma)(t)=0.$$
The above holds provided that
$$(\nabla\DB\cdot uG)(\gamma(t))+
(\nabla\DB\cdot F)(\gamma(t))=0.$$
Assuming that
$$(\nabla\DB\cdot G)(x)\not= 0\qquad\qquad\forall x\in S,\eqno(2.15)$$
the values of the control $u$ are thus uniquely determined by
$$u=\phi(x)\doteq ~ -{\nabla\DB\cdot F(x)\over
\nabla\DB\cdot G(x)}.\eqno(2.16)$$
A turnpike is {\it regular} if for every $x\in S$ (2.15) holds true
and $|\phi(x)|\leq 1$.
A trajectory $\gamma$ is said to be a {\it Z-trajectory}
if there exists a regular turnpike S such that
$\{\gamma(t):t\in Dom(\gamma)\}\subset S$.

An isolated nonordinary arc or INOA $S$ is a {\it barrier} in
$\Omega$ (where $\Omega$ satisfies (C1)-(C4)) if it verifies:
\v
\item {(S1')} For every $x\in S$,
$X(x)$ and $Y(x)$ point to the same side of $S$
\item {(S2')} Each of the function $\Delta_A,\;\Delta_B$ is either
identically zero on $S$ or nowhere zero on $S$.
\v
\n A point $x\in\R^2$ is a {\it near ordinary point} if it is an ordinary
point or belongs to an INOA that is either a turnpike or a barrier.
\v
\n{\bf Remark 2.1.}

The definition of near ordinary point given in Sussmann (1987 a)
is more general, but for our pourpouses this simpler
definition is sufficiently general.

%\null
%\hfuzz=20pt
%\font\bigbf=cmbx10 scaled \magstep3
%\font\medbf=cmbx10 scaled \magstep2
\def \n{\noindent}
\def \v{\vskip 1em}
\def \vs{\vskip 2em}
\def\sgn{\hbox{sgn}}

\def \vsk{\vskip 4em}
\def \M {{\cal M}}

\def\DB{\Delta_B}
\def\DA{\Delta_A}
\def \R {I\!\!R}
\def\ve{\varepsilon}
\def \C {{\cal C}}
\def \P {{\cal P}}

\def\H {{\cal H}}

%\baselineskip = 18 truept
%\input amssym.def
%\input amssym
%\input fonts.tex
%\input macro
%\def\autori{B. PICCOLI}
%\def\titolo{Time-Optimal Control problems for the Swing and the Ski}
%\twelvepoint
\vsk
{\bf 3. Preliminary Theorems.}
\vs
In this section we show three theorems on control systems with control
appearing linearly as in (2.4) and prove a theorem
relating the control systems
for the swing and the ski with this type of control systems.

Consider the control system (2.4). For ordinary points
we have the following:
\v
\n {\bf Theorem 3.1} {\it Let $\Omega\subset\R^2$ be an open set
such that each $x\in \Omega$ is an ordinary point.
Then all time optimal trajectories $\gamma$ for the restriction of (2.4)
to $\Omega$ are bang-bang
with at most one switching. Moreover if $f>0$ throughout $\Omega$ then $\gamma$
is an $X,Y$ or $Y*X-$trajectory, if $f<0$ throughout $\Omega$ then $\gamma$
is an $X,Y$ or $X*Y-$trajectory.}
\v
\n For the proof see Sussmann (1987 a, p.443).
For near ordinary point we have a similar theorem on the local structure
of optimal trajectories, see Sussmann (1987 a, p.459):
\v
\n {\bf Theorem 3.2} {\it Let $x$ be a near ordinary point. Then there
exists a neighborhood $\Omega$ of $x$ such that every time optimal
trajectory $\gamma$ for the restriction of (2.4) to $\Omega$
is concatenation of at most five trajectories each of which
is an $X-,Y-$ or $Z-$trajectory.}
\v
Let consider now a control system as in (2.3):
$$\dot x\ =\ h(x,u).\eqno(3.1)$$
If $w_1,w_2\in\R^2$, define the triangle:
$$C(w_1,w_2)\dot =\{w\in\R^2:w=\lambda w_1+\mu w_2;\lambda,\mu\geq 0;
\lambda+\mu\leq 1\}$$
and consider the condition:
\v
\item {(P1)} There exist two $\C^1$ vector fields
$w^{\pm}(x)$ that for every $x$ either are
linearly independent or have the same versus,
and such that:
$$\{w^{\pm}(x)\}\subset \{h(x,u):u\in U\}$$
$$\{h(x,u):u\in U\}\subset
\{x:x=\lambda w^+(x)+\mu w^{-}(x),\;\lambda,\mu\geq 0,
\;0\leq\lambda+\mu<1\}\cup\{w^{\pm}(x)\}\subset$$
$$\subset C(w^+(x),w^-(x)).$$

\n Suppose that (3.1) verifies (P1). Then we can define a system
of the form (2.4) choosing:
$$F\ =\ {w^++w^-\over 2},\qquad G\ =\ {w^+-w^-\over 2}\eqno(3.2)$$
thus $(F\pm G)(x)=w^{\pm}(x)$.
We can also define the map:
$$\P:\R^2\times U\rightarrow [-1,1]\eqno(3.3)$$
in the following way.
If $w^{\pm}(x)$ are independent then $\P(x,v)=u$
if and only if $h(x,v)$ and $F(x)+uG(x)$ are parallel.
Otherwise, $\P(x,v)$ is constantly equal to $+1$ if $F(x),G(x)$
have the same versus (that is if $|w^+|$ is bigger than $|w^-|$)
and to $-1$ if the opposite happens.
{}From (P1) we have that this map is well defined.

\n Given two points $x,\tilde x\in\R^2$
we consider the two endpoints problem of
steering $x$ to $\tilde x$ in minimum time. Let define the value functions:
$$V(x,\tilde x)\ \dot=\ inf\ \{T(\gamma):
\gamma{\rm\ trajectory\ of\ }(2.4),(3.2), In(\gamma)=x,Term(\gamma)=\tilde
x\}$$
$$\bar V(x,\tilde x)\ \dot=\ inf\ \{T(\gamma):\gamma{\rm\ trajectory\ of\ }
(3.1),In(\gamma)=x,Term(\gamma)=\tilde x\}$$
We have the following:
\v
{\bf Theorem 3.3} {\it Assume that every point is near ordinary for (2.4),(3.2)
and consider two points $x,\tilde x\in\R^2$:
\item {(i)} there exists a trajectory $\gamma$ of (2.4),(3.2)
such that $T(\gamma)=V(x,\tilde x)$

\item {(ii)} if $\Gamma$ is the set of trajectories of (2.4),(3.2)
corresponding to bang-bang controls
and $V'(x,\tilde x)=inf\{T(\gamma):\gamma\in\Gamma,In(\gamma)=x,
Term(\gamma)=\tilde x\}$ then:
$$V'(x,\tilde x)=V(x,\tilde x).$$}
\v
\def \g {\gamma}
{\bf Proof.} Assume that $V(x,\tilde x)<+\infty$, otherwise there is
nothing to prove. The statement $(i)$ follows from the convexity
of the set of velocities for (2.4),(3.2), see for example
Lee and Markus (1967).

\n Let $\gamma$ be a time optimal trajectory steering $x$ to $\tilde x$.
{}From Theorem 3.1 and 3.2 we have that $\gamma$ is concatenation
of $X,Y$ and $Z-$trajectories. If $\gamma$ is bang-bang we are done.
Assume now that $\gamma\restriction [t_0,t_1]$ is a $Z-$trajectory.
{}From the definition of turnpike we have that for every $t\in[t_0,t_1]$
the two vectors $X(\gamma(t)),Y(\gamma(t))$ points to opposite
sides of the image of $\gamma$. If $\gamma(t)$ is sufficiently
near to $\gamma(t_0)$ then we can construct a bang-bang trajectory
$\hat\g$ with one switching steering $\g(t_0)$ to $\g(t_1)$, in the
following way. For $|t-t_0|$ sufficiently small the $Y-$trajectory
$\g^+$ passing through $\g(t_0)$ and the $X-$trajectory $\g^-$ passing through
$\g(t)$ meet each other in at least one point. Let $\bar x$ be the first
point in which $\g^+$ intersect the image of $\g^-$, after passing
through $\g(t_0)$. We can construct $\hat\g$ following $\g^+$ up to the
point $\bar x$ and then $\gamma^-$ up to the point $\g(t)$.
Divide now $[t_0,t_1]$ into $n$ equal subintervals inserting the points
$k_i\dot=t_0+(i/n) (t_1-t_0),i=0,\ldots,n$. If $n$ is sufficiently
large we can construct, as above, a bang-bang trajectory steering
$\g(k_i)$ to $\g(k_{i+1}),i=0,\ldots,n-1$. Let $\hat\g_n$ be the
concatenation of these trajectories. Let $K$ be a compact neighborhood
of the image of $\g$ and let $M\dot=2\;max_K\;(|F|+|G|)<+\infty$.
If $n$ is sufficiently large then by construction $\hat\g_n$ lies
in $K$, hence:
$$|\g(t)-\hat\g_n(t)|\leq {t_0-t_1\over n}\;M\leq
{T(\g)\over n}\;M.$$
The subset $J$ of $Dom(\g)$ on which $\g$ is a $Z-$trajectory is a finite
union of closed interval. Otherwise we can constuct a
sequence $t_n$ in $\partial J$ converging to a time $t\in Dom(\g)$
and then $\g(t)$ is not a near ordinary point.
Repeating the same reasonings for every subinterval of $J$ we obtain
a sequence of bang-bang trajectories converging uniformly to $\g$.
\v
Notice that every bang-bang trajectory for (2.4),(3.2) is a trajectory
for (3.1), then $\bar V(x,\tilde x)\leq V'(x,\tilde x)$.
We can prove the following:
\v
\n{\bf Theorem 3.4} {\it Assume that (P1) holds true and that every
point is a near ordinary point for (2.4),(3.2).
Consider two points $x,\tilde x\in\R^2$ then
\item {(i)} if there exists a bang-bang time optimal
control $u$ for (2.4),(3.2) steering
$x$ to $\tilde x$, then there exists a time optimal control
$v$ for (3.1) corresponding to the same trajectory $\gamma$ of $u$,
i.e. $h(\gamma(t),v(t))=F(\gamma(t))+u(t)G(\gamma(t))\in\{w^{\pm}(\g(t))\}$;

\item {(ii)} if every time optimal control $u$ for (2.4),(3.2)
is not bang-bang (i.e. if $\gamma(t)$ belongs to a turnpike for some interval)
then the time optimal control for (3.1) does not exist but we have:
$$V(x,\tilde x)=\bar V(x,\tilde x)$$
and for each $\ve>0$ there exists a control $v$, corresponding
to a trajectory $\eta$ steering $x$ to $\tilde x$, such that
$h(\eta(t),v(t))\in\{w^{\pm}(\eta(t))\}$ for each $t\in Dom(\gamma)$
and $T(\eta)\leq \bar V(x,\tilde x)+\ve$.

}

\v
{\bf Proof.} Suppose first that there exists a bang-bang time
optimal control $u$
and, by contradiction, that there exists $v$ control of (3.1),
with corresponding trajectory $\eta$ steering $x$ to $\tilde x$, such that
$T(\eta)< T(\gamma)$.
If $h(\eta(t),v(t))\in \{w^{\pm}(\eta(t))\}$
for almost every $t$ then $\eta$
is a trajectory of (2.4),(3.2) contradicting the optimality of $\gamma$.
If this is not true
let define the feedback control $\bar u(\eta(t))\dot=\P(\eta(t),v(t))$,
where $\P$ is the map in (3.3).
There exists a trajectory $\bar\gamma$ corresponding to $\bar u$ that is a
reparametrization of $\eta$ and is a trajectory of (2.4),(3.2).
Finally $|F(\eta(t))+\bar u(\eta(t))G(\eta(t))|\geq
|h(\eta(t),v(t))|$, therefore:
$$T(\eta)\geq T(\bar\gamma)\geq T(\gamma)$$
that gives a contradiction.

Suppose now that every time optimal control $u$ is not bang-bang.
{}From Theorems 3.1 and 3.2 we have that
there exists $I\subset Dom(\gamma)$ such that $\gamma(I)$
is contained in a turnpike. Consider a control $v$ for (3.1)
corresponding to a trajectory $\eta$ steering $x$ to $\tilde x$.
Let define $\bar u,\bar\gamma$ as above
and define $S=\{t:|h(\eta(t),v(t))|
<|F(\eta(t))+\bar u(\eta(t))G(\eta(t))|\}$. If $meas (S)=0$
then $\eta$ is a bang-bang trajectory of (2.4),(3.2)
and is not time optimal.
If $meas(S)>0$ then there exists $n$ and $\delta>0$ such that
$$meas\left(\left\{t:|h(\eta(t),v(t))|
<|F(\eta(t))+\bar u(\eta(t))G(\eta(t))|-{1\over n}\right\}\right)
\geq\delta>0.$$
Therefore:
$$T(\eta)\geq T(\bar\gamma)+\delta\;{1\over n}>T(\gamma)$$
for every time optimal $\gamma$.

\n From Theorem 3.3 it follows
the second part of the statement $(ii)$.
\v
Theorem 3.4 shows the relationships between the two control system
(2.4),(3.2) and (3.1). If we are able to determine the time optimal
control $u$ for a given problem for (2.4),(3.2) then we immediately
know all about the same problem for (3.1). Indeed if the time optimal
control $u$ is bang-bang then the trajectory $\gamma$ of $u$ is
a trajectory also for (3.1), is optimal and corresponds to a control
$v$ taking values in $\{r_-,r_+\}$.
If every time optimal control is not bang-bang
then the time optimal control for (3.1) does not exists.

\n Thank to Theorems 3.1 and 3.2 in most cases we are able to know
the local structure of time optimal control for (2.4),(3.2) and
then for (3.1). Using the same methods of Sussmann (1987 a,b,c)
and Piccoli (1993 a,b), we are
able to solve explicitely many optimization problems for the swing
and the ski models. Some examples will be given in the following sections.

\null
\hfuzz=20pt

\font\medbf=cmbx10 scaled \magstep2
\def \n{\noindent}
\def \v{\vskip 1em}
\def \vs{\vskip 2em}
\def\sgn{\hbox{sgn}}
\def \vsk{\vskip 4em}

\def\DB{\Delta_B}
\def\DA{\Delta_A}
\def \R {I\!\!R}
\def\ve{\varepsilon}
\def \C {{\cal C}}
\def \M {{\cal M}}
\def \P {{\cal P}}

%\input amssym.def
%\input amssym
%\input fonts.tex
%\input macro
%\def\autori{B. PICCOLI}
%\def\titolo{Time-Optimal Control Problems for the Swing and the Ski}
%\twelvepoint
%\vfill\eject
\vsk
\n{\bf 4.Time Optimal Control of the Swing.}
\vs
In this section we treat the problem of time optimal control for
the swing considering the minimum time problem and a Mayer type problem.

Recall equation (1.2).
We have to verify that this control system satisfies the condition (P1)
of section 3. If $sin(x_1)=0$ then $\dot x_2=0$ for every $v\in [r_-,r_+]$
and the set of velocities is a segment. If $sin(x_1)\not= 0$,
from the second equation of (1.2) we obtain:
$$v\ =\ -{\dot x_2\over g\;sin(x_1)}\eqno(4.1)$$
and replacing (4.1) in the first equation of (1.2):
$$\dot x_1\ =\ {x_2\;g^2\;sin^2(x_1)\over\dot x_2^2}.\eqno(4.2)$$
It is easy to check from (4.2) that the set of velocities for (1.2)
at a given point lies on a branch of hyperbola. Defining
$w^{\pm}=h(x,r_{\mp})$, (P1) holds. Notice that with this definition
we have:
$$(F+G)(x)=h(x,r_-)\qquad (F-G)(x)=h(x,r_+).$$
Therefore if $u$ is a bang-bang control for the auxiliary system
(2.4),(3.2) with corresponding trajectory $\g$, then
$\g$ is also a trajectory for (3.1) and corresponds to the control:
$$v(t)=r_-\quad {\rm if}\quad u(t)=1,\qquad
  v(t)=r_+\quad {\rm if}\quad u(t)=-1.$$

Hence we can compute the linear system (3.2) associated to this control system
as in section 3, explicitely:
$$F=\left(\matrix{{r_+^2+r_-^2\over 2\;r_+^2\;r_-^2}\;x_2 &\cr
                    -{g\over 2}\;(r_++r_-)\;sin(x_1)        &\cr}\right)
\qquad
G=\left(\matrix{{r_+^2-r_-^2\over 2\;r_+^2\;r_-^2}\;x_2 &\cr
                    {g\over 2}\;(r_+-r_-)\;sin(x_1)        &\cr}\right).$$
For simplicity we define:
$$a\;\dot=\;(r_+\;r_-)^2\quad b\;\dot=\;r_++r_-\quad c\;\dot=\;r_+-r_-
\quad d\;\dot=\;r_+^2+r_-^2\eqno(4.3)$$
then:
$$F=\left(\matrix{{d\over 2\;a}\;x_2 &\cr
                    -{g\over 2}\;b\;sin(x_1)        &\cr}\right)
\qquad
G=\left(\matrix{{b\;c\over 2\;a}\;x_2 &\cr
                    {g\over 2}\;c\;sin(x_1)        &\cr}\right)
.\eqno(4.4)$$
To investigate the local structure of time optimal trajectories
we have to compute the functions $\DA,\DB$, defined in (2.10,11):
$$\DA={g\;c\;(d+b^2)\over 4\;a}\;x_2\;sin(x_1)\eqno(4.5)$$
$$[F,G]=\left(\matrix{-{g\;c\;(d+b^2)\over 4\;a}\;sin(x_1) &\cr
                    {g\;c\;(d+b^2)\over 4\;a}\;b\;x_2\;cos(x_1)&\cr}\right)$$
hence:
$$\DB={g\;c^2\;(d+b^2)\over 8\;a}\ \left({b\over a}\;x_2^2\;cos(x_1)+g\;
sin^2(x_1)\right).\eqno(4.6)$$
Every turnpike is subset of $\DB^{-1}(0)$, then we have to solve
the equation $\DB(x)=0$ that is equivalent to:
$${b\over a}\;x_2^2\;cos(x_1)+g\;sin^2(x_1)\ =\ 0$$
that gives:
$$x_2^2\ =\ -{g\;a\over b}\;{sin^2(x_1)\over cos(x_1)}.$$
By periodicity we can restrict ourselves to the case $x_1\in[0,2\pi]$.
There is the isolated solution $(0,0)$ and
if $x_1\in]{\pi\over 2},{3\pi\over 2}[$ we have the two solutions:
$$x_2\ =\ \pm\sqrt{-{g\;a\over b}\;{sin^2(x_1)\over cos(x_1)}}\eqno(4.7)$$
otherwise there is no solution.

\n The two branches of solutions form two curves, one contained
in the first quadrant, the other in the forth quadrant. The two
curves meet each other at the point $(\pi,0)$ and $(0,0),(\pi,0)$
are the only points of $\DB^{-1}(0)$ that verify $\nabla\DB(x)=0$.
For each $x\in\DB^{-1}(0)\setminus \{(0,0),(\pi,0)\}$
we can compute the control $\phi(x)$ defined in (2.8):
$$\phi(x)\ =\ \phi(x_1)=
\ {(2\;b^2-d)\;cos^2(x_1)-d\over b\;c\;(3\;cos^2(x_1)+1)}.
\eqno(4.8)$$
{}From (4.8) we have:
$$\lim_{x_1\to {\pi\over 2}}\phi(x_1)\ =\
\lim_{x_1\to {3\pi\over 2}}\phi(x_1)\ =\
-{r_+^2+r_-^2\over r_+^2-r_-^2}\ <\ -1$$
$$\lim_{x_1\to \pi}\phi(x_1)\ =\ {2\;b^2-2\;d\over 4\;b\;c}\ =
\ {r_+\;r_-\over r_+^2-r_-^2}>0$$
$${d\phi(x_1)\over dx_1}\ =\ -{4\;(b^2+d)\;sin(x_1)\;cos(x_1)\over
b\;c\;(3\;cos^2(x_1)+1)^2}$$
then $\phi$ is increasing for $x_1\in ]\pi/2,\pi[$ and decresing
for $x_1\in ]\pi,3\;\pi/2[$. Therefore there exist $\ve_1>0,\ve_2\geq0$
such that $|\phi(x_1)|\leq 1$ for $x_1\in ]\pi/2+\ve_1,\pi-\ve_2[
\;\cup\; ]\pi+\ve_2,3\;\pi/2-\ve_1[\not=\emptyset$.
We have that $\ve_2=0$ if and only if:
$${r_+\;r_-\over r_+^2-r_-^2}\leq 1$$
i.e. if and only if:
$$r_+\geq r_-\;{1+\sqrt{5}\over 2}.$$
Hence regular turnpikes do exist. We can choose two points
$x,\tilde x\in\R^2$ such that the only time optimal
trajectory for (2.4),(4.4) that
steers $x$ to $\tilde x$ is not bang-bang: it is sufficient to take
two points of the same turnpike, see Sussmann (1987 a).
In this case the time optimal
control for (1.2) does not exist and the second part of Theorem 3.4
applies.
\v
Now, Theorems 3.1 and 3.2 determine the local structure of the time
optimal trajectories for the system (2.4),(4.4) and then for the
swing model.

\n Observe that $\DA^{-1}(0)=\{(x_1,x_2):x_2=0{\rm\ or\ }
\sin (x_1)=0\}$. Therefore it easy to check that $\DA^{-1}(0)
\backslash\{(0,0),(\pi,0)\}$ is union of a finite number of
INOAs that are barriers. Moreover, from the reasoning above we have
that $\DB^{-1}(0)\backslash\{(0,0),(\pi,0)\}$ is union of a finite
number of INOAs each of which is either a turnpike or a barrier.
Hence every point of $\Omega\dot=\R^2\backslash\{(0,0),(\pi,0)\}$
is near ordinary. Notice that if $\g$ is a trajectory of (2.4),(4.4)
and $\g(t)\in\{(0,0),(\pi,0)\}$ for some time $t\in Dom(\g)$, then
$\g$ is a constant trajectory. Indeed, in this case, the control
has no effect being $X(\g(t))=Y(\g(t))=0$.
Moreover, $(0,0),(\pi,0)$ are the only point in which either $X$
or $Y$ vanishes. Thus if $\g$ is a trajectory of (2.4),(4.4) we have
that either $\g$ is constant or $\g(t)\in\Omega$ for every
$t\in Dom(\g)$. In the latter case we can apply Theorems 3.1 and 3.2.

\n Given the expression of $\DA,\DB$ we can calculate $f$ of (2.14):
$$f(x_1,x_2)\ \dot=\ {c\left[{b\over a}\,x_2^2\,\cos(x_1)+g\,\sin^2(x_1)
\right]\over 2\,x_2\,\sin(x_1)}.$$
Consider the set:
$$Q\dot= \left\{(x_1,x_2):|x_1|\leq{\pi\over 2}\right\}$$
and a time optimal trajectory $\g$
that verifies $\g(t)\in Q$ for every
$t\in Dom(\g)$. We have, from Theorems 3.1 and 3.2, that $\g$ is
bang--bang and that $\g$ can switches from control $+1$ to
control $-1$ if $x_2\,\sin(x_1)>0$ and from control $-1$ to control $+1$
if $x_2\,\sin(x_1)<0$. Therefore $Q\backslash\DA^{-1}(0)$
is divided is four parts and on each parts only one kind of
switching is permitted. This correspond to the fact that
if the swing is raising his distance from the earth then the swinger
can change only from control $r_-$ to control $r_+$, instead if the swing
is lowering his distance from the earth then the swinger can change
only from $r_+$ to $r_-$.

\n In this case, the map $\P$ of (3.3) is bijective for
almost every $x$ and we can
establish a bijective correspondence between the trajectories of (1.2) and
of (2.4),(4.4). In particular if $\g:[a,b]\rightarrow \R^2$
is a trajectory of (2.4),(4.4) then there exists a trajectory
$\eta:[c,d]\rightarrow\R^2$ of (1.2) verifying $\eta(c)=\g(a)$,
$\eta(d)=\g(b)$ and $\eta(t)\in \g([a,b])$ for every $t\in [c,d]$.
Therefore also some geometric properties of (2.4),(4.4) hold for (1.2).
For example the reachable sets from a given point are the same
for (2.4),(4.4) and for (1.2).

\n In Alberto Bressan (1993)
it was considered the problem fo raising the amplitude of the
first half oscillation starting from a given point $x\in\R^2$.
Notice that if $\sin(x_1)>0$ then $F_2-G_2<F_2+G_2<0$
and the opposite happens if $\sin(x_1)<0$.
Hence, comparing the vector $X,Y$, it is easily seen that if we want to reach
the maximum amplitude of the first half oscillation we must
choose the control $+1$ if $x_2\,\sin(x_1)>0$ and $-1$ if
$x_2\,\sin(x_1)<0$. For the swing this means to choose
$r_-$ if we are raising the distance from the earth
and $r_+$ in the other case.
Obviously this is also the optimal control to raise the amplitude
after a given number of oscillation.

Using the Pontryagin Maximum Principle and Lemma 2.1,2.2 and 2.3 we can prove
the following:
\v
\n {\bf Lemma 4.1.} {\it Assume that $\eta:[a,b]\rightarrow\R^2$
is a bang--bang trajectory of (1.2),
$\eta(t)\in Q$ for every $t\in [a,b]$,
and that there exists $t_1,t_2\in (a,b)$ such that
either the first or the second component of $\eta(t_1)$ vanishes
and the same hold for $\eta(t_2)$.
If $\eta$ has no switching then $\eta$ can not be optimal.}
\v
{\bf Proof.} Assume for example that $\eta$ corresponds to the constant
control $+1$ and, by contradiction, that $\eta$ is optimal.
Observe that $\eta$ is a trajectory of (2.4),(4.4), then
if $\eta$ is optimal there exists an adjoint covector field
$\lambda(t),\ t\in [a,b]$, along $\eta$. We have, from Lemma 2.2, that
$\lambda(t_1)\cdot G(\eta(t_1))\geq 0$ and,
from $\DA(\eta(t_1))=0$, that $G(\eta(t_1)),Y(\eta(t_1))$ are parallel.
{}From Lemma 2.3 and (4.6) the function:
$$\alpha(t)\dot= arg\Big(Y\big(\eta(t_1)\big),\
v\big(G(\eta(t)),t;t_1\big)\Big)$$
is strictly increasing. It is easy to check
from (2.5) and (4.4) that the function:
$$\psi(t)\dot=det\Big(v\big(G(\eta(t_2)),t_2;t\big),\
v\big(Y(\eta(t_2)),t_2;t\big)\Big)$$
is constant. Assume for example that $t_2>t_1$. We have that
$G(t_2)=v(G(\eta(t_2)),t_2;t_2)$
and $Y(t_2)=v(Y(\eta(t_2)),t_2;t_2)$
are parallel because $\DA(\eta(t_2))=0$. Hence
$v\big(G(\eta(t_2)),t_2;t_1\big)$ and $v\big(Y(\eta(t_2)),t_2;t_1\big)$
are parallel, but from (2.5) it follows $v\big(Y(\eta(t_2)),t_2;t_1\big)
=Y(\eta(t_1))$. We conclude that $\alpha(t_2)=k\,\pi$ for
some integer $k>0$.

\n Now, from Lemma 2.1 and 2.2:
$$\lambda(t)\cdot G(\eta(t))=
\lambda(t_1)\cdot v\big(G(\eta(t)),t;t_1\big)\geq 0\quad\forall t\in [a,b].$$
But since $\alpha(t_2)\geq\pi$ the vector $v\big(G(\eta(t)),t;t_1\big)$
for $t\in [a,b]$ makes a rotation of an angle strictly greater than $\pi$
and then there is no vector $\lambda(t_1)$ for which the above
inequality can hold.
\v
It is useful the following:
\v
\n{\bf Lemma 4.2.} {\it Assume that $\g$ is a time optimal trajectory
of (2.4),(4.4), that $\g$ has a switching
at time $t_1\in Dom(\g)$ and that $\DA(\g(t_1))=0$.
Then $\DA(\g(t_2))=0$, $t_2\in Dom(\g)$, if and only if $t_2$ is
a switching time for $\g$.}
\v
{\bf Proof.} Let $u(t),\ t\in Dom(\g)$, be the control corresponding
to $\g$. Since $\g$ is optimal there exists an adjoint covector
field $\lambda$ along $(u,\g)$.
Let $t_2$ be the first time that either is a switching time or that verifies
$\DA(\g(t_2))=0$. We have that $u\restriction [t_1,t_2]$ is constant,
say $u\equiv 1$, and that $\lambda(t_1)\cdot G(\g(t_1))=0$.
But $G(\g(t_1))$ and $F(\g(t_1))+G(\g(t_1))$
are parallel, hence $\lambda(t_1)\cdot [F(\g(t_1))+G(\g(t_1))]=0$.
{}From Lemma 2.1 have:
$$\lambda(t_2)\cdot \big[F(\g(t_2))+G(\g(t_2))\big]=
\lambda(t_1)\cdot v\big(F(\g(t_2))+G(\g(t_2)),t_2;t_1\big)=$$
$$=\lambda(t_1)\cdot \big[F(\g(t_1))+G(\g(t_1))\big]=0.$$
Now if $\DA(\g(t_2))=0$ we have that
$G(\g(t_2))$,$F(\g(t_2))+G(\g(t_2))$ are parallel then
$t_2$ is a switching time. On the other hand, if $t_2$ is a switching time
then $\lambda(t_2)\cdot G(\g(t_2))=0$, hence being $\lambda(t_2)\not=0$
we have that $G(\g(t_2))$ and $F(\g(t_2))+G(\g(t_2))$ are parallel.
This means $\DA(\g(t_2))=0$. We can argue in the same way for the other
switching times.
\v
Consider now the problem of reaching, with the swing, an angle
$\bar \theta\in(0,\pi/2)$ in minimum time with given initial condition.
We argue about the corresponding problem for (2.4),(4.4)
with initial point $\bar x\in Q$. We can solve this problem constructing
all the time optimal trajectories starting from $\bar x$ and
using the (final) transversality condition of the PMP,
see for example Lee and Markus (1967),
to select among these trajectories the optimal one.
To construct all the time optimal trajectories we can proceed as in
Piccoli (1993 a,b).

\n Assume, for example, that $\bar x=(\bar x_1,\bar x_2)$,
$\bar x_1>0,\bar x_2>0$, being similar the other cases.
Let $\g$ be a time optimal trajectory that verifies $In(\g)=\g(0)=\bar x$.
If $\g$ lies on the first quadrant for $t\in[0,t_1]$ then
$\g\restriction [0,t_1]$ is a $Y-$trajectory or a $X*Y-$trajectory.
Let $\g^+$ be the $Y-$trajectory satisfying $\g^+(0)=\bar x$,
and let $t^+$ be the first positive time, if any, for which $\g^+(t^+)$
belongs to the $x_1-$axis.

\n First we want to prove that if $t^+<\infty$,
$\g(0)=\bar x$ and $\g\restriction
[0,t^++\ve]$ is a $Y-$trajectory for some $\ve>0$, then $\g$
is not time optimal after a given time.
Assume, by contradiction, that $\g$ is time optimal.
{}From Theorem 3.1, being $f>0$ in the fourth quadrant, we have that
necessarily $\g$ is a $Y-$trajectory until it reaches the $x_2-$axis,
say at time $t_1$. Now $\g$ can not switch at time $t_1$ because
of Lemma 4.2. Hence we can apply Lemma 4.1 obtaining a contradiction.

\n Then if we want to contruct all the time optimal trajectories
it is enough to consider the trajectories $\g_s$, $\g_s(0)=\bar x$,
that follow $\g^+$ for a given time $s\in[0,t^+]$ and then switch
to control $-1$. Now if $\g_s$ is time optimal then there exists
an adjoint covector field $\lambda_s$ along $\g_s$.
Since $\lambda_s(s)\cdot G(\g_s(s))=0$ we can determine $\lambda_s$
up to the product by a positive scalar.
Hence we can compute the switching function
$\phi_s=\lambda_s\cdot G(\g_s))$ and, by Lemma 2.2, determine
the behaviour of $\g_s$, that is its switching times.
Using Lemma 4.2 we can see that, for every $s$,
$\g_s$ has to switch on each quadrant.
More precisely $\g_s$ will make a second switching before reaching
the $x_2-$axis, then a third switching before the second time of
intersection with the $x_1-$axis and so on. The set of switching
points of $\g_s$ form some manifolds called
switching curves, see Piccoli (1993 b) for the exact definition. Each switching
curve lyes on a quadrant. After a given number of oscillations
some $\g_s$ reach the manifold $\M \dot=
\{(x_1,x_2):x_1=\bar\theta\}$.
To select the optimal trajectory between these $\g_s$, we can use the
transversality condition that, in this case, is $\lambda_s(t_s)
\cdot (0,1)=0$ if $t_s$ is the first time of intersection of $\g_s$
with $\M$.

\n To solve numerically our problem, first we have to solve the
equation for the swing (1.2), for constant control $r_+$ and $r_-$,
and a class of initial data (that is to approximate some elliptic integrals).
Indeed every $\g_s$ is a finite concatenation of such trajectories.
Then we have to consider the complementary equation (2.6) for $\g_s$.
This is a linear system and can be solved numerically by usual methods,
determining, up to a scalar,
the adjoint covector field $\lambda_s$ with initial condition
$\lambda_s(s)\cdot G(\g_s(s))=0$. Finally the transversality condition
determines a value $\bar s$,
and then the corresponding trajectory $\g_{\bar s}$.

\n If $\bar x\notin Q$ and $|\bar\theta|>\pi/2$ then the construction
of time optimal trajectories is more difficult. We can do it
following Piccoli (1993 a,b),
but in this case we have to take into account
the turnpikes. It can happen that some time optimal trajectories are
$Z-$trajectories for some time interval of positive measure, and the optimal
control for the swing does not exist. An example of this situation
is given in the following.

\v
We now consider the Mayer problem:
$$\cases{\dot \g(t)\;=\;F(\g(t))+u(t)\;G(\g(t))   &\cr
         \g(0)=\bar x               &\cr
         {\rm min\ }\{t:\exists\gamma{\rm\ such\ that\ }
                 \gamma(t)\in {\M}\}  &\cr}\eqno(4.9)$$
where ${\M}\dot =\{(x_1,x_2):x_1=c\},c\in\R$. We consider the case
in which $c=3\;\pi/2$, $\bar x=(\bar x_1,\bar x_2)$ belongs to a turnpike and
$\bar x_1\in ]\pi,3\pi/2[,\bar x_2>0$. Our aim is to
show the existence of some values of $r_{\pm}$ such that
the problem (4.9) with the dynamics (1.2) does not have a solution.

If $\gamma$ is a solution then in particular $\gamma$ is time
optimal for steering $\bar x$ to $Term(\gamma)$ and hence is a concatenation
of $X-,Y-$ and $Z-$trajectories.
Following the algorithm in Piccoli (1993 a,b),
we can cover a region of the plane
with the time optimal trajectories starting from $\bar x$.
Let $S$ be the turnpike
to which $\bar x$ belongs. The region of first quadrant below $S$ is
covered by $Y-$trajectories originating from $S$ and
the region over $S$ is covered by
$X-$trajectories. If $x'$ is the endpoint of $S$ that comes after
$\bar x$ for the orientation given by $X,Y$, then the $X-$trajectories
that cross $\DB^{-1}(0)$ over $x'$ must switch changing the control
to $+1$. The points in which these trajectories change control
form a curve called switching curve, see Piccoli (1993 b).
Therefore if $\gamma$ is a time optimal trajectory we have
two possibilities: either $\gamma=\gamma_2*\gamma_1$ where $\gamma_1$ is
a $Z$-trajectory possibly trivial and $\gamma_2$ is a $Y-$trajectory,
or $\gamma=\gamma_3*\gamma_2*\gamma_1$ where $\gamma_1$ is
a $Z$-trajectory possibly trivial, $\gamma_2$ is an $X-$trajectory
and $\gamma_3$ is a $Y-$trajectory.

Our aim is to find some $r_{\pm}$ such that the solution to (4.9) is
of the first type with $\gamma_1$ not trivial, i.e. $Dom(\gamma)$
is not a single point.

Let define $\gamma_x^+$ as the $Y-$trajectory that verifies
$In(\gamma_x^+)=\gamma_x^+(0)=x$ and let $\gamma_x^S$ be
the $Z-$trajectory that satisfies $In(\gamma_x^S)=\gamma_x^S(0)=x$.
We have that $\gamma_x^+$ satisfies:
$$\cases{\dot x_1={d+b\;c\over 2\;a}\;x_2   &\cr
         \dot x_2={g\;(c-b)\over 2}\;sin(x_1) &\cr}$$
hence:
$$\ddot x_1=\left({d+b\;c\over 2\;a}\right)\;\dot x_2=
\left({d+b\;c\over 2\;a}\right)\;\left({g\;(c-b)\over 2}\right)\;
sin(x_1).\eqno(4.10)$$
We define:
$$\alpha\ \dot=\ {d+b\;c\over 2\;a},\qquad
\beta\ \dot=\ {g\;(c-b)\over 2}\eqno(4.11)$$
and $\omega$ by:
$$\omega^2\ \dot=\ -\alpha\;\beta.\eqno(4.12)$$
We can solve (4.10) in $t$ using the first integral:
$$\dot x_1^2-2\;\omega^2\;cos(x_1)$$
and obtaining:
$$\dot x_1^2=\alpha^2\;\bar x_2^2+4\;\omega^2\;\left(
sin^2\;{\bar x_1\over 2}-sin^2\;{x_1\over 2}\right).\eqno(4.13)$$
Now let $\tilde x\dot=\gamma_{\bar x}^+(T)$ where $T$ is the first
time at which $\gamma_{\bar x}^+$ intersects $\M$ and define:
$$k\ \dot=\ sin\;{\bar x_1\over 2}\qquad a_0^2\ \dot=\
{4\;\omega^2\over\alpha^2\;\bar x_2^2+4\;\omega^2\;k^2}.\eqno(4.14)$$
After straightforward calculations from (4.13,14) we obtain:
$$T\ =\ {1\over\sqrt{\alpha^2\;\bar x_2^2+4\;\omega^2\;k^2}}\
\int_{\bar x_1}^{\tilde x_1}\;{dx_1\over\sqrt{1-a_0^2\;sin^2\;
{x_1\over 2}}}$$
and, using the substitution $\theta\dot={x_1\over 2}$:
$$T\;=\;{2\over\sqrt{\alpha^2\;\bar x_2^2+4\;\omega^2\;k^2}}\;
\int_{{\bar x_1\over 2}}^{{\tilde x_1\over 2}}
\;{d\theta\over\sqrt{1-a_0^2\;sin^2\;\theta}}\;=$$
$$=\;{2\over\sqrt{\alpha^2\;\bar x_2^2+4\;\omega^2\;k^2}}\;
\left[E\left(a_0,{\tilde x_1\over 2}\right)-E\left(a_0,{\bar x_1\over 2}
\right)\right]\eqno(4.15)$$
where:
$$E(l,\bar\theta)\ =\ \int_{0}^{\bar\theta}\;
{d\theta\over\sqrt{1-l^2\;sin^2\;\theta}}\eqno(4.16)$$
is the elliptic integral of the second type.
Now we want to compute the time along the trajectory $\gamma_s$
that satisfies $\gamma_s\restriction [0,s]=\gamma_{\bar x}^S\restriction
[0,s]$ and that is a $Y-$trajectory after the time $s$.
Define $T(s)$ to be the first time in which $\gamma_s$ intersects
${\M}$. We have:
$$T(s)\ =\ s+{2\over\sqrt{\alpha^2\;x_2^2(s)+4\;\omega^2\;k_s^2}}\;
\left[E\left(a_s,{3\;\pi\over 2}\right)-E\left(a_s,{x_1(s)\over 2}
\right)\right]\eqno(4.17)$$
where:
$$x(s)\ \dot=\ \gamma_{\bar x}^S(s)\quad k_s\ \dot=\ sin\left(
{x_1(s)\over 2}\right)\quad a_s^2\ \dot=\
{4\;\omega^2\over\alpha^2\;x_2^2(s)+4\;\omega^2\;k_s^2}.\eqno(4.18)$$
For any $s\geq 0$ we can calculate the difference between
$T$ and $T(s)$. It is clear that if:
$$\left.{d\over ds}T(s)\right|_{s=0}\eqno(4.19)$$
is negative then $\gamma_{\bar x}^+$ is not optimal for (4.9).
In computing (4.19) we will use the approximations:
$$x_1(s)=\bar x_1+(F_1+\phi G_1)(\bar x)\;s+o(s)\qquad
x_2(s)=\bar x_2+(F_2+\phi G_2)(\bar x)\;s+o(s).$$

After some calculations from (4.15,16,17,18) it follows:
$$\left.{d\over ds}\;T(s)\right|_{s=0}=$$
$$1+\left(\alpha^2\;\bar x_2^2+4\;\omega^2\;k^2\right)^{-{3\over 2}}\;A\;
\left[-\int_{\bar x_1\over 2}^{\tilde x_1\over 2}
\left(1-a_0^2\;sin^2\;\theta\right)^{-{1\over 2}}+
\left(1-a_0^2\;sin^2\;\theta\right)^{-{3\over 2}}\;a_0^2\;sin^2\;\theta\;
d\theta\right]+$$
$$-{d+\phi(\bar x_1)\;b\;c\over 2\;a}\ \bar x_2\
\left[\left(\alpha^2\;\bar x_2^2+4\;\omega^2\;k^2\right)\;
\left(1-a_0^2\;k^2\right)\right]^{-{1\over 2}}\eqno(4.20)$$
where:
$$A\;=\;\alpha^2\;g\;\left(\phi(\bar x_1)\;c-b\right)\;\bar x_2\;sin\;\bar x_1
+2\;\omega^2\;\bar x_2\;sin\;{\bar x_1\over 2}\;cos\;{\bar x_1\over 2}\;
\left({d+\phi(\bar x_1)\;b\;c\over a}\right)=$$
$$=\;{\alpha\;g\;c\over 2\;a}\;\left(d+b^2\right)\;
\left(\phi(\bar x_1)-1\right)\;\bar x_2\;sin\;\bar x_1.\eqno(4.21)$$

Notice that the second and the third terms on the righthand side
of (4.20) are negative. If $r_-$ tends to zero then
the second term tends to:
$$-\ {cos(\bar x_1)\over sin(\bar x_1)}\
{(3\,cos^2(\bar x_1)-cos(\bar x_1)+2)\over 3\,cos^2(\bar x_1)+1}\
\left({\tilde x_1-\bar x_1\over 2}\right).\eqno(4.22)$$
Now if we choose $r_+$ sufficiently large then $|\phi(x_1)|\leq 1$
for $x_1\in ]\pi,\pi+\ve[$ and some $\ve>0$, hence we can take
$\bar x_1$ arbitrarily near to $\pi$. But as $\bar x_1$ tends
to $\pi$, the expression in (4.22) tends to minus infinity.
In particular we can choose $r_-,r_+$ and $\bar x_1$ in such a way that:
$$\left.{d\over ds}T(s)\right|_{s=0}<0.\eqno(4.23)$$
Now it is easy to see that, choosing $r_{\pm},\bar x_1$ as above,
the $Y*X-$trajectory leaving from $\bar x$
can not be optimal for (4.9), in fact $X(\bar x)$ tends to zero.
On the other hand for (4.23) the $Y-$trajectory starting
from $\bar x$ can not be optimal.
Hence the time optimal trajectory for (4.9) contains
a non trivial $Z-$trajectory. Using Theorem 3.4 we obtain that the
problem (4.9) with the dynamics (1.2) has no optimal solution.

Given the time optimal control $u$ for (4.9), corresponding
to the trajectory $\gamma$, we can consider the feedback control
$v$ such that, see (3.3):
$$\P\big(\gamma(t),v(\gamma(t))\big)\ =\ u(t).$$
Indeed, in this case, the function ${\cal P}$, defined in the third section,
is bijective outside the coordinate axes and the line $\{x_1=\pi\}$.
We have:
$$v(\gamma(t))\ =\ \root 3\of{-a\;{u(t)\;c-b\over u(t)\;b\;c+d}}.\eqno(4.24)$$
The trajectory $\eta$ corresponding to the control $v$ is a
reparametrization of $\gamma$.

\n Since $\dot\gamma_1(t),\dot\eta_1(t)>0$
for every $t$ we can consider the inverse functions $t_{\gamma}(x_1),
t_{\eta}(x_1)$ and define $\gamma(x_1)=\gamma(t_{\gamma}(x_1))$, $\eta(x_1)=
\eta(t_{\eta}(x_1))$. If $\tilde x=Term(\gamma)$ then:
$$T(\gamma)\ =\ \int_{\bar x_1}^{\tilde x_1}\;{dx_1\over
F_1\big(\gamma(x_1)\big)+u(t_{\gamma}(x_1))\;G_1\big(\gamma(x_1)\big)}$$
$$T(\eta)\ =\ \int_{\bar x_1}^{\tilde x_1}\;{v^2(\eta(x_1))dx_1\over
\eta_2(x_1)}.$$
Then we can calculate the difference of times $T(\gamma)-T(\eta)$
computing the difference of the two integrals. We can obatin
a better performance using a bang-bang control but it is not easy
to compute explicitly one such control and
its total variation tends to infinity as its time tends to the minimum.
The control $v$ has the advantage of being defined by the explicit formula
(4.24) and of having a fixed variation.

\null
\hfuzz=20pt

\font\medbf=cmbx10 scaled \magstep2
\def \n{\noindent}
\def \v{\vskip 1em}
\def \vs{\vskip 2em}
\def\sgn{\hbox{sgn}}

\def \vsk{\vskip 4em}
\def \M {{\cal M}}

\def\DB{\Delta_B}
\def\DA{\Delta_A}
\def \R {I\!\!R}
\def\ve{\varepsilon}
\def \C {{\cal C}}

\def\H {{\cal H}}

%\input amssym.def
%\input amssym
%\input fonts.tex
%\def\autori{B. PICCOLI}
%\def\titolo{Time-Optimal Control Problems for the Swing and the Ski}
%\baselineskip = 18 truept

\vsk
%\twelvepoint
%\input macro
\def \I{{\cal I}}
\n {\bf 5.Time Optimal Control for the  Ski.}
\vs
In this section we deal with time optimal control problems for the ski.
We consider a skier
on a one-dimensional trail described by a curve $(x(s),y(s))$ in the plane.
The skier can choose the height $u$ of his baricentre, that
is the distance from the point of contact with the trail,
within certain physical bounds.
We refer to Aldo Bressan (1991) for assumptions and notations.
The dynamics is given by
the system:
$$\cases{\dot s\ =\ {p\over {\I} (s,u)}     &\cr
                          & $u\in [r_-,r_+]$ \cr
         \dot p\ =\ {{\I}_s\over 2\;{\I}^2}\;p^2-M\;g\;(1-c\;u)\;y'(s)&\cr}
         \eqno(5.1)$$
where $s$ is the arclength of the trail, $p$ its conjugate momentum,
$c(s)$ the curvature of the trail,
$\I$ is $c^2$ times the inertial moment of the pair ski-skier with respect
to the centre of curvature,
${\I}_s$ its partial derivative with respect to $s$,
$g$ the gravity acceleration and
$M$ is the total mass of the pair ski-skier.
We approximate the body of the skier with a system of
$n$ equal masses with height $i\;h/n,\;i=1,\ldots,n,$ where $h$
is the height of the top of the skier. We have:
$$h=b\,u,\qquad b=2\;{m+m_s\over m}\;{n\over n+1}$$
and, see Bressan (1991):
$$\I (s,u)\ =\ \alpha\;c^2(s)+m_s+{m\over n}\;
\sum_{i=1}^n\;\left(1-c(s)\;b\;u\;{i\over n}\right)^2\eqno(5.2)$$
where $\alpha$ is the inertial moment of the ski
with respect to its baricentre, $m_s$ its mass and
$m$ the mass of the skier, so that $M=m+m_s$. In Bressan (1991)
it was assumed that:
$$u\;c\ \leq\ 1\qquad u\in [r_-,r_+]\eqno(5.3)$$
$$sgn({\partial \I\over\partial u})=-sgn(c)\eqno(5.4)$$
where, by definition, $sgn(0)=0$.
It is easy to verify that (5.3,4) hold under the hypotheses (5.2).

We now consider the case in which the curvature $c$ is constant.
Then ${\I}_s=0$ and (5.1) reduces to:
$$\cases{\dot s\ =\ p\;\left(\alpha\;c^2(s)+m_s+m+m\;c^2\;b^2\;u^2\;
{(2n+1)(n+1)\over 6n^2}-m\;c\;b\;u\;{(n+1)\over n}\right)^{-1}   &\cr
\dot p\ =\ -M\;g\;(1-c\;u)\;y'(s)                  &\cr}.\eqno(5.5)$$
Now we want to verify that, under suitable conditions on the curvature,
the control system (5.5) satisfies the hypothesis (P1) of the third section.
If $c=0$ then the control disappear and we obtain a dynamical system.
If $p=0$ or $y'(s)=0$ the set of velocities is a segment and we are done.
Hence from now on we assume that $c,p,y'(s)\not=0$.
{}From the second equation of (5.5) we can express $u$ as a function
of $\dot p$ and substitute its value in the first equation obtaining:
$$\dot s\ =\ {p\over P(\dot p)}$$
where:
$$P(\dot p)=\dot p^2\;\left({m\,b^2\,(2n+1)\,(n+1)
\over M^2\,g^2\,[y'(s)]^2\,6n^2}\right)
+\dot p\;\left({m\,b\,[b(2n+1)-3n]\,(n+1)
\over M\,g\,y'(s)\,3n^2}\right)$$
$$+{6\,\alpha\,c^2\,n^2+m\,b^2\,(2n^2+3n+1)-6\,m\,b\,n\,(n+1)+6\,(m+m_s)\,n^2
\over 6\,n^2}.$$
The discriminant of the polynomial $P$ is:
$$disc(P)=\;-\ {m\,b^2\,(n+1)\,[2\,\alpha\,c^2\,(2n+1)+m\,(n+1)
+2\,m_s\,(2n+1)]
\over 3\,M^2\,g^2\,[y'(s)]^2\,n^2}$$
then $P(\dot p)$ has no zeros.

\n The zeros of the second derivative of $\dot s$ with
respect to $\dot p$ are:
$$\dot p_{\pm}={g\,M\,y'(s)\,[3n-b\,(2n+1)]\over b\,(2n+1)}\;\hfill$$
$$\hfill
\pm\;{g\,M\,y'(s)\,n\over \sqrt{m}\,b\,(2n+1)\,(n+1)}\;
\sqrt{2\,(\alpha\,c^2+m_s)\,(2\,n^2+3\,n+1)+m\,
(n^2-1)}.$$
Then we compute $u_{\pm}\in\R$ in such a way that $\dot p(u_{\pm})=
\dot p_{\pm}$:
$$u_{\pm}\;=\;{n\,\left[3\,\sqrt{m}\,(n+1)
\pm\,\sqrt{2\,(\alpha\,c^2+m_s)\,(2\,n^2+3\,n+1)+m\,(n^2-1)}\right]
\over b\,c\,\sqrt{m}\,(2n+1)\,(n+1)}.$$
If $c<0$ then $u_+<0$ and the condition $u_-<0$ gives:
$$c^2\ <\ {m\;(4\;n+5)\over \alpha\;(2\;n+1)}\;-\;{m_s\over\alpha}.
\eqno(5.6)$$
If $c>0$ then (5.6) gives $u_->0$ but if we assume:
$${3\,n\over b\,c\,(2n+1)}-r_+>0\eqno(5.7)$$
then the condition $u_-> r_+$ gives:
$$a_2\;c^2+a_1\;c+a_0>0\eqno(5.8)$$
where:
$$a_2=(2n+1)\,(n+1)\,[m\,b^2\,r_+\,(2n+1)\,(n+1)-2\,\alpha\,n^2]\eqno(5.9)$$
$$a_1=-6\,m\,b\,r_+\,n\,(2n+1)\,(n+1)^2
\quad a_0=2\,n^2\,(n+1)\,[m\,(4n+5)-m_s\,(2n+1)].$$
For $n$ sufficiently large and for standard data,
the conditions (5.7,8,9) are similar
to the condition (5.3) only slightly more restrictive.
If we assume that (5.6,7,8,9) are verified then
$\dot s(\dot p(u)), u\in [r_+,r_-]$ is convex or concave, depending
on the signes of $p,c$.
Moreover we have that $\dot s(\dot p)$ tends to zero as $|\dot p|$
tends to infinity and it is easy to check that (P1) is verified.

The assumptions (5.6,7,8,9) are natural in fact in the application
we have some bounds on the curvature.
\v
We can now compute the associated system with control appearing
linearly, as in section 3, to obtain
information about the time optimal controls. If $c=0$ the
control does not appear and we have a dynamical system as observed above.
Hence from now on we assume that $c\not= 0$. We have:
$$F=\left(\matrix{{p\;({\I}_++{\I}_-)\over 2\;{\I}_+\;{\I}_-}\cr
                 -{M\;g\over 2}\;y'(s)\;(2-c\;(r_++r_-)) \cr}\right)\qquad
  G=\left(\matrix{{p\;({\I}_+-{\I}_-)\over 2\;{\I}_+\;{\I}_-}\cr
                 -{M\;g\over 2}\;y'(s)\;c\;(r_+-r_-) \cr}\right)\eqno(5.10)$$
where ${\I}_{\pm}={\I}(r_{\pm})$. After straightforward calculations
it follows:
$$[F,G]\ =\ {m\;g\over 2}\;L\;\left(\matrix{y'(s)\cr -p\;y''(s)\cr}\right)$$
where:
$$H\dot= {{\I}_++{\I}_-\over 2\;{\I}_+\;{\I}_-}\quad
J\dot= {{\I}_+-{\I}_-\over 2\;{\I}_+\;{\I}_-}\quad
L\dot=c\;H\;(r_+-r_-)-J\;(2-c\;(r_++r_-))$$
and:
$${\Delta}_B\ =\ {m\;g\over 4}\;\big(-2\;J\;L\;p^2\;y''(s)
+m\;g\;L\;c\;(y'(s))^2\;(r_+-r_-)\big).$$
Therefore the equation for turnpikes is:
$$p^2\ =\ {m\;g\;c\;(r_+-r_-)\over 2\;J}\;{(y'(s))^2\over y''(s)}.
\eqno(5.11)$$
Observe that the first factor of (5.11) is negative in fact $sgn(J)=-sgn(c)$
(from (5.4)). For example if we consider the curve:
$$\left({cos(|c|\,s)\over |c|},{sin(|c|\,s)\over |c|}\right),
\quad s\in [0,2\pi/|c|]\eqno(5.12)$$
we have:
$$p\ =\ \pm\sqrt{{m\;g\;(r_+-r_-)\over 2\;|J|}\;{cos^2(|c|\,s)\over
sin(|c|\,s)}}.$$
We can argue as in section 4 for time optimal control problems.
Given a minimum time problem,
if every time optimal trajectory for the system (2.4),(5.10) contains
a $Z-$trajectory then the time optimal control for the ski does not
exist. If the opposite happens then the time optimal control exists and
is bang-bang.
The geometric structure of the system (2.4),(5.10) is very similar to that
one of the system (2.4),(4.4). The only difference is the sign
of the function $f$ of (2.14) that can be different.
Consider again the curve (5.12).
Since we are considering a skier, (5.12) makes sense
only if we consider the restriction $s\in [0,\pi/|c|]$ if $c<0$
and $s\in[\pi/|c|,2\,\pi/|c|]$ if $c>0$.

\n For the case $c>0$ we can repeat the reasoning made for the swing,
that is for the system (2.4),(4.4). Indeed, with the above restriction,
the system (2.4),(5.10) is similar to the system (2.4),(4.4) restricted
to the set $Q$ of section 4.
In particular we have
that the time optimal control to reach a given position $\bar s$,
with fixed initial data, in minimum time
does exist and is bang--bang. The natural hypothesis is that
the skier has at least one choice that permit him to reach $\bar s$
without turning back, this means $\dot s(t)>0$ for every $t$.
Hence we have that the optimal trajectory lies in the part of the plane
$\{(s,p):p>0\}$. Therefore the time optimal control has at most two switchings
and can be determined as in section 4.

\n If $c<0$ we are in the opposite case because the natural
restriction force us to stay in the part of the plane in which
there are turnpikes. The time optimal control does not necessarily exist.
Again we can repeat the reasoning of section 4.

\n For Mayer problem (4.9) we have obtained the nonexistence assuming
some conditions on $r_{\pm}$, but it can happen that these conditions
are not phisically acceptable for the ski model.

\vfill
\eject
%\bye

\null
\hfuzz=20pt

\font\medbf=cmbx10 scaled \magstep2
\font\ninesc=cmcsc10 at 12pt
\def \n{\noindent}
\def \v{\vskip 1em}
\def \vs{\vskip 2em}
\def \vsk{\vskip 4em}

\def \R {I\!\!R}
\def\ve{\varepsilon}
%\baselineskip = 18 truept
%\input amssym.def
%\input amssym
%\input fonts.tex
%\input macro
%\def\autori{B. PICCOLI}
%\def\titolo{Time-Optimal Control problems for the Swing and the Ski}
%\twelvepoint
\centerline{\medbf References}
\vs
\n\item {} {\ninesc Bressan, Alberto}, 1993, Impulsive Control Systems,
to appear on {\it Proceedings of IMA Workshop on Geometrical Control Theory}.
\v
\item {} {\ninesc Bressan, Aldo}, 1991, On Some Control Problems concerning
the Ski or Swing, {\it Atti Accad. Naz. Lincei, Memoirs}, {\bf IX-1},
149--196.
\v
\item {} {\ninesc Bressan, Aldo, and Motta, M.}, 1994, On Control Problems
of Minimum Time for Lagrangean Systems similar to a Swing, to appear.
\v
\item {} {\ninesc Lee, E.B., and Markus, L.}, 1967,
{\it Foundations of Optimal Control Theory}, (New York:Wiley).
\v
\item {} {\ninesc Piccoli, B.}, 1993 a, An algorithm for the Planar Time
Optimal
Feedback Synthesis, Preprint SISSA-ISAS;
1993 b, Classification of Generic Singularities
for the Planar Time Optimal Synthesis, Preprint SISSA-ISAS.
\v
\item {} {\ninesc Sussmann, H.J.}, 1987 a, The Structure of Time--Optimal
Trajectories for Single--Input Systems in the Plane: the
$\C^{\infty}$ Nonsingular Case,
{\it SIAM J. Control and Opt.}, {\bf 25}, 433--465;
1987 b, The Structure of Time--Optimal
Trajectories for Single--Input Systems in the Plane: the
General Real Analytic case,
{\it SIAM J. Control and Opt.}, {\bf 25}, 868--904;
1987 c, Regular Synthesis for Time--Optimal
Control of Single--Input Real--Analytic Systems in the Plane,
{\it SIAM J. Control and Opt.}, {\bf 25}, 1145--1162.

\bye